\documentclass[12pt]{article}

\usepackage{amsmath}
\usepackage{mathtools}
\usepackage{xcolor}
\usepackage{amssymb}
\usepackage{amsfonts}
\usepackage{graphicx}
\usepackage{booktabs}
\usepackage{units}
\usepackage{xcolor}
\usepackage{comment}
\usepackage{caption}
\usepackage{subcaption} 
\usepackage{dblfloatfix}   
\usepackage{algorithm}
\usepackage{algorithmicx}
\usepackage{algpseudocode}
\algrenewcommand\algorithmicindent{0.5em}

\usepackage{hyperref}
\hypersetup{
    colorlinks=true,
    linkcolor=blue,     
    urlcolor=blue,
}
\usepackage{xparse}

\usepackage{tikz}
\usepackage{pgfplots}
\usepackage{appendix}

\usepackage{enumitem}

\newcommand{\R}{\mathbb{R}}

\floatname{algorithm}{Pseudo code}

\newtheorem{definition}{Definition}

\newcounter{thmcounter}[section]  
\renewcommand{\thethmcounter}{\thesection.\arabic{thmcounter}}

\newcommand{\defthmwithqed}[2]{%
  \NewDocumentEnvironment{#1}{ o }{%
    \refstepcounter{thmcounter}
    \begin{trivlist}%
      \item[\hskip \labelsep \bfseries #2~\thethmcounter%
        \IfValueT{##1}{\ (##1)}.]%
  }{%
    \hfill $\square$%
    \end{trivlist}%
  }%
}

\defthmwithqed{theorem}{Theorem}
\defthmwithqed{lemma}{Lemma}
\defthmwithqed{assumption}{Assumption}
\defthmwithqed{modelass}{Modeling Assumption}
\defthmwithqed{defi}{Definition}
\defthmwithqed{prop}{Proposition}
\defthmwithqed{prob}{Problem}
\defthmwithqed{cor}{Corollary}
\defthmwithqed{remark}{Remark}

\title{Neural Scaling Laws for Learning-based Identification of Nonlinear Systems}
\author{Hannes Gernandt, Marco Roschkowski, Karim Cherifi \\
\normalsize University of Wuppertal \\
\normalsize \texttt{\{gernandt, roschkowski, cherifi\}@uni-wuppertal.de}}

\date{\today}

\begin{document}

\maketitle
\tableofcontents

\begin{abstract}
The use of machine learning models in system identification has increased due to their ability to approximate complex nonlinear dynamics with high accuracy. However, often it is not clear how the performance of trained models scales with given resources such as data, compute, and model size. To allow for a better understanding of the scalability of the performance of machine learning models, we verify neural scaling laws (NSLs) in the context of system identification from input-state-output data using different evaluation metrics for accuracy and different system architectures, including input-affine and physics-informed port-Hamiltonian representations. Our verified NSLs can help to forecast performance improvements and guide model design or data acquisition.
\end{abstract}

\textbf{Keywords:} System identification, neural scaling laws,  port-Hamiltonian systems, nonlinear systems, control applications

\maketitle

\section{Introduction} \label{sec:intro}
System identification is a discipline in systems and control theory that deals with the building of  mathematical models of dynamical systems from data \cite{Ljung1999}. {While methods relying on linear models and carefully designed experiments provide strong guarantees and interpretability \cite{Ljung1999}, they often fall short when dealing with the complex, nonlinear dynamics that appear in many modern engineering and scientific applications.}

In recent years, the increasing availability of data and computational power has opened the door to machine learning-based approaches to system identification. Neural networks and other {highly expressive} function approximators have shown great promise in capturing complex nonlinear dynamics that are difficult to model with classical methods. {Nonlinear system identification is particularly important for modern applications \cite{Nelles2020} such as robotics~\cite{Av21,Cat26}, autonomous driving~\cite{Ro21,Nie22,Hua20}, and climate modeling \cite{Ou24}, where dynamics are inherently complex and linear approximations are often insufficient}. Alongside purely data-driven approaches, there has also been a surge of interest in physics-informed learning methods, which aim to integrate domain knowledge and physical priors into the identification process, thereby improving sample efficiency, robustness and interpretability.

Despite these advances, fundamental challenges remain in the identification of nonlinear systems. Training large {neural networks based models} is computationally expensive, data collection may be costly or time-consuming, and there are no clear rules for determining whether a given identification task is feasible under resource constraints. These bottlenecks motivate the need for neural scaling laws (NSLs), which are {empirical} predictive rules that relate the accuracy of the model to the amount of data, computation, and model capacity. In machine learning, NSLs have emerged as a powerful tool to guide practitioners in allocating resources for training large models \cite{hestness2017deep}. 
{More precisely, for a given machine learning task we associate an NSL as a function $L:R \longrightarrow [0, \infty)$ with the interpretation that $L(r)$ provides a lower bound for the achievable evaluation error by a trained model using resources $r \in R$. The resources $r$ can be measured in terms of compute, dataset size or model capacity, leading to a distinct neural scaling law in each case.}

NSLs were first studied in 
\cite{hestness2017deep} with the key finding that these {often} take the very simple form of a saturated power law \begin{equation}\label{eq:power_law}
    L(r) = \alpha + \beta r^\delta
\end{equation}
with $\alpha, \beta \ge0, \delta < 0$.
In \cite{hestness2017deep}, it was also pointed out that a certain amount of minimal resources must be committed to observe this scaling behavior. With a very small amount of resources {such as the amount of data and model size}, the model will only be able to learn ``best guessing".
On the $\log$-$\log$ scale, \eqref{eq:power_law} is the graph of a function that is approximately linear with slope $\delta$ representing smooth scaling. The scaling then saturates once $r$ is large enough such that $\alpha \sim \beta r^\delta$.

NSLs have been observed for several machine learning tasks such as language models~\cite{kaplan2020scaling} where NSLs have been used to predict the resources necessary to successfully train a $175$~billion parameter large language model~\cite{brown2020language}.
In addition, NSLs were determined for image-text models in \cite{alabdulmohsin2022revisiting}, for vision transformers in  \cite{zhai2022scaling}, and for language image pretraining \cite{cherti2023reproducible}. Recently, feature learning has been investigated as a way to improve NSLs in~\cite{bordelon2025feature}.
There have also been theoretical efforts in understanding NSLs \cite{maloney2022solvable,bahri2024explaining,paquette20244+}. {Since these theoretical investigations are not yet general enough to cover complex machine learning systems typical in learning-based system identification, we will focus on a purely empirical approach.}

In this paper, we bring these ideas to the domain of learning-based system identification. Specifically, we present the first systematic experimental study of NSLs in this context. 
{We study several state-space system architectures like input-affine and port-Hamiltonian for which the coefficients that are represented by multi-layer perceptrons (MLPs) as well as input-output approaches like NARX and SUBNET.} Our results show how scaling behavior manifests itself in system identification tasks and we provide practical information on the resources required to achieve the desired accuracy. We believe that these findings offer a principled way to evaluate the feasibility of identification problems and to plan resource allocation.

By training thousands of system identification models, we obtain a large set of outcomes consisting of the evaluation metrics and data size, compute and model capacity for each of those trained models. This set of outcomes allows us to determine the NSLs of the corresponding combination of dataset and identification method. 
We observe that NSLs in system identification may have several breaks as laid out by \cite{caballero2022broken} and, therefore, we {use their definition} by considering  NSLs of the form 
\begin{align}
\label{brokenNSL_intro}    
    L(r) = \alpha + \beta \ r^{\delta_0} \prod\limits_{i = 1}^b \left[ 1 + \left(\frac{r}{\sigma_i}\right)^{\frac{1}{\varphi_i}}\right]^{\delta_i \cdot \varphi_i}
\end{align}
for a suitable choice of real parameters $\alpha,\beta,\sigma_i,\varphi_i,\delta_i$ and the number of breaks $b\geq 0$. Note that for $b=0$ breaks we recover the saturated power law~\eqref{eq:power_law}.

Our contributions and findings can be summarized as follows:
\begin{itemize}
\item 
{We {show that broken NSLs}~\eqref{brokenNSL_intro} {are present} in the context of learning based system identification and {empirically} verify them for various examples and system architectures.}

\item {Our results indicate that early trends in the scaling curve could allow us to anticipate the data, compute, and model complexity required for higher accuracy based on the scaling law. This can help to forecast performance improvements and guide model design or data acquisition in learning-based system identification methods.}
 \item {Our results indicate that data availability is more critical than model size. Increasing the number of parameters does not necessarily lead to improved performance, as the NSL values frequently saturate at relatively small model sizes. This behavior contrasts with observations in the natural language processing literature, where performance often scales with model size \cite{brown2020language}.}

  \item Training of system models for an unexpectedly large number of epochs can still reduce the evaluation error.
  \item {The NSL results are largely independent of the chosen model structure. Overall, we observe consistent behavior across models with different levels of incorporated physical priors.}
\end{itemize}

The rest of the paper is organized as follows. In Section~\ref{sec:SysID}, we recall learning-based system identification methods and system architectures that are used in this paper. In Section~\ref{sec:NSL}, we define a general class of NSLs that we use in the context of system identification. In Section~\ref{sec:approximation} we present our experimental setup to generate the set of outcomes and describe how the NSLs are approximated from these outcomes. In Sections~\ref{sec:verification_ios} and~\ref{sec:verification_io} , we verify compute, data, and model NSLs for several example systems. 

\section{Learning-based system identification} \label{sec:SysID}

In order to study NSLs in system identification,
we consider a family of models that
can approximate a broad range of control systems while still allowing for systematic variation in data, model size, and training compute.  

\subsection{Identification from input-state-output data}
Specifically, we consider three architectures of {state space models of} control systems ranging from unstructured models to physics-based modeling approach:
\begin{itemize}
    \item \textit{unstructured} systems without structural restrictions \cite{Nelles2020,CHEN90,NarP90}
    \begin{align}\tag{3a}
    \label{eq:unstructured}
    \dot{x}(t) = f(x(t), u(t)),\quad  
    y(t)= g(x(t)),
    \end{align}
with state $x(t)\in\mathbb{R}^n$, input $u(t)\in\mathbb{R}^m$, output  $y(t)\in\mathbb{R}^l$ and 
for some $f:\R^n\times\R^m\rightarrow\R^n$ and $g:\R^n\rightarrow\R^l$;
    \item \textit{input-affine} systems \cite{Xie25}
    \begin{align}\tag{3b}
        \label{eq:input_affine}
    \dot{x}(t) = h(x(t)) + j(x(t)) u(t), \quad  y(t) = k(x(t)),
    \end{align}
    where $h:\R^n\rightarrow\R^n$, $k:\R^n\rightarrow\R^l$ and  $j:\R^n\rightarrow\R^{n\times m}$;
    \item \textit{port-Hamiltonian} (pH) systems \cite{cherifi2025}
    \begin{align}
    \nonumber
        \dot{x}(t) &= (J(x(t)) - R(x(t))) \nabla H(x(t)) + B(x(t)) u(t), \\
    y(t) &= B(x(t))^\top \nabla H(x(t)),   \tag{3c} \label{eq:ph}
    \end{align}
    with matrix-valued functions $J,R:\R^n\rightarrow\R^{n\times n}$, $B:\R^n\rightarrow\R^{n\times m}$, a~continuously differentiable Hamiltonian $H:\mathbb{R}^n\rightarrow[0,\infty)$, with gradient $\nabla H(x)\in\mathbb{R}^n$, and with $J(x)=-J(x)^\top$ and $R(x) = R(x)^\top \geq 0$.
\end{itemize}
Compared to \eqref{eq:unstructured} and \eqref{eq:input_affine}, the structural constraints on the system coefficients in the pH~representation~\eqref{eq:ph} guarantee that the learned model preserves physical properties such as passivity, energy dissipation or energy conservation \cite{SchJ14}.

Our approach to learning-based system identification is based on neural network parameterizations of dynamical systems where we represent all unknown functions in \eqref{eq:unstructured}, \eqref{eq:input_affine} and \eqref{eq:ph} by multilayer perceptrons (MLPs) to enable the approximation of general nonlinear dependencies. 
In pH~systems, maintaining the pH structure of the functions in~\eqref{eq:ph} requires customized parameterizations, as outlined by~\cite{cherifi2025}.

In the learning-based system identification, we assume that there is a~sufficiently rich collection of input-state-output data that is generated based on a {suitable sampling of initial values and excitation signals  \cite{cherifi2025,Neary2023,vanOMWTJS24} which is described in more detail in Section~\ref{sec:approximation}. In particular, it is assumed that the state-space dimension and the values of the derivatives are known exactly and that there is no noise.}

For the data sets given, the functions of all three architectures (unstructured, input-affine, pH) given by \eqref{eq:unstructured}, \eqref{eq:input_affine}, and \eqref{eq:ph} are approximated by MLPs as explained in~\cite{cherifi2025}. 
For the training of the MLPs, we use the same MSE-type loss function as in \cite{cherifi2025}, that is normalized over individual batches and we normalize the states and inputs before processing those with the neural networks.

In the following, we collect all generated trajectories for all initial values $x_0$ and all input functions $u$ and over a certain time horizon with time-step size $\Delta t$ in matrices {discretized states $x_i\in\mathbb{R}^n$ at time $i\Delta t$ and we consider the entries of $x_i=(x_{ij})_{j=1}^n$ which are referred to as the $j$th state at time $i\Delta t$. The right hand side of the state-space dynamics is then used to obtain the matrix of derivatives of the states} $\dot{X}=(\dot x_{ij}) \in \mathbb{R}^{K \times n}$ where $K \in \mathbb{N}$ denotes the number of data points. The resulting trajectories that are obtained after training for the same initial values and input functions are denoted by $\widehat{\dot{X}}=(\hat{\dot{x}}_{ij}) \in \mathbb{R}^{K \times n}$. {As evaluation metrics, we consider the normalized mean absolute error ($\operatorname{nMAE}$) and the normalized mean squared error ($\operatorname{nMSE}$) defined by the formulas} \begin{align}
\label{eq:score_MAE}
    \operatorname{nMAE}(\dot{X}, \widehat{\dot{X}}) &= \frac{1}{K{\cdot n}} \sum\limits_{i=1}^K \sum\limits_{j = 1}^n \frac{|\dot{x}_{ij} - \widehat{\dot{x}}_{ij}|}{\sigma_j},\\
\label{eq:score_MSE}
    \operatorname{nMSE}(\dot{X}, \widehat{\dot{X}}) &= \frac{1}{K{\cdot n}} \sum\limits_{i=1}^K \sum\limits_{j = 1}^n \left(\frac{\dot{x}_{ij} - \widehat{\dot{x}}_{ij}}{\sigma_j}\right)^2,
\end{align}
where $\sigma_j$ is the standard deviation of the derivative of the {$j$th entry of the state vector} on the entire dataset
{
\[
\sigma_j =  \sqrt{ \frac{1}{K-1}  \sum\limits_{i=1}^K \biggl(\dot{x}_{ij} - \frac{1}{K} \sum\limits_{k=1}^K \dot{x}_{kj} \biggr)^2 }.
\]
If we instead plug in $1$ for $\sigma_j$ for all $j=1,\ldots,n$, then we recover the \textit{mean absolute error} ($\operatorname{MAE}$) and the \textit{mean squared error} ($\operatorname{MSE}$). }
 Note that normalization {with $\sigma_j$} is necessary because a plain $\operatorname{MSE}$ or $\operatorname{MAE}$ metric would not equally balance the error on all states if these have different magnitudes.

\subsection{Identification from input-output data}
{Since state data is not available in most applications, we also consider two methods purely based on input-output data which is given as a finite sequence $((u_1,y_1),\ldots,(u_N,y_N))$ in $\mathbb{R}^m\times\mathbb{R}^l$.} 
\subsubsection{Nonlinear Autoregressive with Exogenous Inputs (NARX)}
{The first architecture, we consider NARX models which are a classical method to predict the output in the next time-step $y_t$ from a lagged window of sizes $n_y,n_u>0$ consisting of past outputs $(y_{t-1}, \dots, y_{t-n_y})$ and past inputs $(u_{t}, \dots, u_{t-n_u})$, see \cite{Nelles2020,CHEN90}. This means that the NARX architecture is of the form
\[
y_t = F(y_{t-1}, \dots, y_{t-n_y}, u_t, u_{t-1}, \dots, u_{t-n_u})
\]
where we model $F$ via an MLP, since NSLs are a phenomenon specifically related neural networks. 
\subsubsection{SUBNET}
NARX does not produce a state-space representation of the system, which might lead to less accurate models. Therefore, we consider the recent SUBNET architecture, see \cite{beintema2023deep} which considers an additional latent state dynamics. 
We briefly explain the architecture of SUBNET, for a more detailed explanation, we refer to~\cite{beintema2023deep}. 
In SUBNET we consider a small input-output trajectory of length $T \in \mathbb{N}$ that is generated by randomly choosing $t \in \{1, \dots, N-T+1\}$ and selecting
$u_t, u_{t+1}, \dots u_{t+T-1}$ as the sequence of inputs and $y_t, y_{t+1}, \dots y_{t+T-1}$ as the sequence of outputs used to train the latent dynamics that aims to forecast step-by-step
\begin{align}
    \label{eq:latent}
\hat{x}_{t+k, t} = f(\hat{x}_{t+k-1, t}),\quad k \in \{1, \dots, T-1\}
\end{align}
for some neural network $f$. In addition, SUBNET, uses a \textit{subspace encoder} which is another neural network that generates  an initial condition $\hat{x}_{t, t}$ for the latent dynamics~\ref{eq:latent}. Finally, a third neural network $h$ is used to predict the corresponding output values
$\hat{y}_{t+k, t} = h(\hat{x}_{t+k, t})$. Note that the particular selection procedure of $t$ and $T$ accelerates the training process since $T \ll N$ and still allows to use most of the training data at each training step if the batch size is large enough. 
}

{
Besides this basic version of SUBNET which uses discrete-time data, there is also a continuous-time version called CT-SUBNET~\cite{sztaki10596} and a more structured port-Hamiltonian version called output-error port-Hamiltonian neural network (OE-PHNN)~\cite{moradi2026port}.}

{
As evaluation metrics, we consider the $\operatorname{nMAE}$ and $\operatorname{nMSE}$ on the output, which is necessary since state data is not available in this case. More precisely, assume that $u_1, \dots, u_{N_t}$ is the sequence of inputs in the test data and $y_1, \dots, y_{N_t}$ the corresponding outputs. Assume that $y_{ij}$, $i \in \{1, \dots, N_t\}$, $j \in \{1, \dots, l\}$ is the sequence of outputs on the test dataset where $l$ is the number of outputs and $N_t$ the length of the testing trajectory.
Then set $Y=(y_{ij}) \in \mathbb{R}^{N_t\times l}$ and \begin{align*}
    \operatorname{nMAE}(Y, \hat{Y}) &= \frac{1}{N_t \cdot l} \sum\limits_{i=1}^{N_t} \sum\limits_{j=1}^{l} \frac{|y_{ij} - \hat{y}_{ij}|}{\sigma_j} \\
    \operatorname{nMSE}(Y, \hat{Y}) &= \frac{1}{N_t \cdot l} \sum\limits_{i=1}^{N_t} \sum\limits_{j=1}^{l} \biggl(\frac{y_{ij} - \hat{y}_{ij}}{\sigma_j}\biggr)^2
\end{align*}
with $\hat{Y} \in \mathbb{R}^{N_t \times l}$ being the predictions of the model.
}

\section{Neural scaling laws gfor system identification} \label{sec:NSL}
NSLs quantify the intuitive dependence between the model validation error given by the $\operatorname{nMAE}$ metric \eqref{eq:score_MAE} or the $\operatorname{nMSE}$ metric \eqref{eq:score_MSE} and the used \textit{resources} $r>0$, where $r$ represents either the \textit{size of the dataset} $d$, the \textit{model parameters} $p$, or the \textit{compute used during training}~$c$ in $\operatorname{flops}$.

We will now give a definition of NSL that is tailored to our application to system identification. To derive NSL we select one of the resources $r$ (data, compute or model size) and vary over this resource and the remaining training parameters. 
For the resource $r$ and each such parameter configuration, we perform a training run for the model and compute  the validation error $e$ of the learned model based on \eqref{eq:score_MAE} or~\eqref{eq:score_MSE}. {Repeating the training runs for each resource $r_i$ and recording the corresponding validation errors $e_i$ produces} to a set of \textit{outcomes} $O=\{(r_1,e_1),(r_2,e_2),\ldots,(r_J,e_J)\}$ in~$\mathbb{R}_+^2$.

For the approximation of NSL, we assume that without restriction the entries of $O$ are ordered with respect to the size of the resources, i.e.\ $r_1\leq r_2\leq\ldots\leq r_J$ and we define the \textit{lower envelope} ${E_O}: [r_1, r_J] \to[0,\infty)$ given by 
\begin{align}
\label{envelope}
{E_O}(r) := \min \{e_i~|~i=1, \ldots, J,\, r_i \le r,\, {(e_i,r_i)\in O}\}. 
\end{align}
{which constitutes a lower bound on the observed errors $e_i$ among all pairs of  outcomes $(r_i,e_i)$ that require not more than $r$ resources.} 

{We define the lower envelope $E_O$ as a function since we aim to interpolate $E_O$ in the following.} To this end, we generate additional interpolation points by choosing a grid size $K \in \mathbb{N}$ and 
\begin{align}
    \label{interpolation_points}
{r_1} = {\tilde{r}_1}  \leq \tilde r_2 \leq  \dots \leq  \tilde r_K = r_J,\quad \tilde r_k := {r_1}^{1 - \frac{k}{K}} r_J^{\frac{k}{K}}.
\end{align}
Since the resources $r_i$ typically vary at exponential scale, the choice of the interpolation points $\{\tilde r_i\}_{i=1}^K$ achieves an equidistant partitioning on $\log$-scale. Furthermore, for each interpolation point, {we set $\tilde e_k =E_O(\tilde{r}_k)$. This leads to the set $\{(\tilde r_i,\tilde e_i)\}_{i=1}^K$ of interpolation points to fit the lower envelope $E_O$.}
\begin{definition}\label{def: broken nsl}
Consider a set of outcomes $O\subseteq\mathbb{R}_+^2$, the lower envelope $E_O$ defined as~\eqref{envelope} and interpolation points $\{\tilde r_i\}_{i=1}^K$ given by~\eqref{interpolation_points}.
Then $O$ obeys a \textit{NSL with at most $b$ breaks and margin $M>0$} if there are $\alpha, \beta \ge 0 \ge \delta_i $ and $ \sigma_i,\varphi_i > 0$, $b \in \mathbb{N}_0$ such that the function $L: [r_1, r_J] \to \mathbb{R}$ given by
\begin{equation}
\label{eq:broken_NSL}
    L(r) = \alpha + \beta \ r^{\delta_0} \prod\limits_{i = 1}^b \left[ 1 + \left(\frac{r}{\sigma_i}\right)^{\frac{1}{\varphi_i}}\right]^{\delta_i \cdot \varphi_i}
 \end{equation}
satisfies \begin{equation}\label{eqn:margin}
    M = \frac{1}{K}\sum\limits_{k=1}^K \left( \log \hat{e}_k - \log \tilde{e}_k \right)^2 
\end{equation}  with $\hat{e}_k = L(\tilde r_k)$ and $\tilde{e}_k={E_O}(\tilde{r}_k)$.
\end{definition}
Compared to existing definitions for NSLs in the literature, we use a broken NSL~\cite{caballero2022broken} and quantify in Definition~\ref{def: broken nsl} an approximation margin $M$ of the NSL which is  particularly motivated by the applications to system identification. In this context, NSL have to be verified for different architectures \eqref{eq:unstructured}, \eqref{eq:input_affine}, \eqref{eq:ph},  various datasets for different application examples and the different resource types data, compute, and model parameters. Here, Definition~\ref{def: broken nsl} allows us to (semi-)automatically compute the NSL and validate its margin $M$ in a well-defined manner independent of the particular setup.
Note that other choices of the lower envelope ${E_O}$ are also possible. For example, by interpolating between the nearest points in $O$. However, for our obtained set of outcomes this does not make a huge difference since there is typically only a small gap between $e_i$ and $e_{i+1}$. We also point out that equally spacing the $e_i$ at $\log$-scale is necessary for \eqref{eqn:margin} to be meaningful. Indeed if we would have simply used all the points in $O$, the margin would be highly biased towards regions containing a high density of configurations.

\section{Generation of outcomes and approximation of neural scaling laws}\label{sec:approximation}
In this section, we explain how the sets of outcomes are generated and how the NSLs were fitted to the observed data relating resources and validation error of a large number of training runs.

\subsection{Model training using input-state-output data}\label{subsec:model_training}
We train the models on the same nonlinear example systems as in \cite{cherifi2025}:
\begin{itemize}
    \item A mass-spring-damper chain, called  \textit{spring system}, which contains nonlinear dissipation;
    \item a magnetically levitated ball system, called \textit{ball system}, which contains nonlinear dissipation and {total energy};
    \item a permanent magnet synchronous motor (PMSM), called  \textit{motor system}, which contains nonlinear {coupling terms}.
\end{itemize}
A detailed description of these systems is given in the appendix.
\subsubsection{Data generation}
{For the generation of trajectory data, we uniformly sample initial states $x_{i,0} \sim U(x_{i, \operatorname{min}}, x_{i, \operatorname{max}})$, $i=1,\ldots,n$, where the choice of $x_{i, \operatorname{min}}$ and $x_{i, \operatorname{max}}$ depends on the particular application, and we consider input signals {to excite the system}
    \begin{align}
        \label{def:input}
    u(t) =
    \sum_{k=1}^N
    a \sin(2\pi k f_0t +\phi_k), 
    \end{align}
    for suitable parameters $f_0,N,a>0$  and uniformly random sampled phases $\phi_k\in[0,2\pi)$. The sampled trajectories are then generated using 
SciPy’s \texttt{odeint} function. We assume that we know the governing state-space equations to generate the data and hence, in particular, the state space dimension and that there is no noise during the data generation.}
\subsubsection{Training setup}
We vary the number of trajectories, training epochs, depth of the model and hidden dimensions of the model and use the same training script to gather the training evaluations necessary to determine NSLs.
The number of epochs used to train the model runs over a grid of $n_e = 2, 4, \dots, 2^{14}$.

One full trajectory represents a $10$ seconds long simulation of the system with time discretization step $\Delta t = 0.01$ seconds which implies that one trajectory consists of $1000$ data points.
We first randomly select $\tilde{d} = 2, 4, \dots, 2^9$. We then select a random float $d \in [\tilde{d}/2, \tilde{d}]$. The number of training trajectories is then $n_t = \frac{d}{10}$. Since $n_t$ does not need to be an integer, we also consider fractional trajectories that are computed by $n_t = 10 \cdot a + b$, $a \in \mathbb{N}, b \in [0, 1)$. We then generate $a+1$ trajectories and use the first $a$ full trajectories. Finally, we select the smallest integer $\tilde{b}$ which is smaller than $b \cdot 1000$ and use the first $\tilde{b}$ data points from the last trajectory.
It is important to consider fractional trajectories since the interval $[2^0, 2^1]$ contains very few integers. Hence, considering fractional trajectories allows one to explore significantly more configurations for small data sizes.

The hidden dimension of the model is obtained by running a loop with values $\tilde{n}_h=2, 4, 8, 16$ and selecting a random integer $n_h$ between $\tilde{n}_h/2$ and $\tilde{n}_h$. The reason is again that we want to consider as many model sizes as possible.
Similarly, the depth of the model is obtained by running a loop for the depth with values $\tilde{n}_d=2, 4$ and selecting a random integer $n_d$ between $\tilde{n}_d/2$ and $\tilde{n}_d$. As in~\cite{cherifi2025}, the hidden dimension for the MLP is multiplied by~$2$ and rounded for unstructured models, to compensate for structured models using more sub-modules, roughly equalizing the parameter count at a given hidden dimension. The total parameter count then equals approximately $4\cdot n_d \cdot n_h^2$.
Finally, all these loops are nested in the lines 4-7 of Algorithm~\ref{algo:submit}, resulting in the submission of $9\cdot 14\cdot 4 \cdot 2 = 1008$ experiments being performed per dataset/architecture combination. We find that there is a failure rate of approximately $3\%$ {due to freezes in the program execution}. Thus, the outcomes of such experiments are not considered{, because the final error scores are not available in these cases}.

We use random seeds that determine the random selection of ${d}, n_h$ and $n_d$ to ensure reproducibility of our experiments and vary the seed to obtain more experiments with different random choices and we use $N=7$ different random seeds to obtain a sufficiently large set of outcomes $O$.
{In line 12 of Algorithm}~\ref{algo:submit}, the expression $\operatorname{main}(n_e, n_t, n_d, n_h, A, D))$ is an abbreviation for submitting a job where a model with depth $n_d$, hidden dimension $n_h$ is trained with the number of training trajectories $n_t$ and number of training epochs $n_e$ on the dataset $D$ with architecture $A$. This function then returns the metrics $\operatorname{nMAE}, \operatorname{nMSE}$ and the number of parameters of the model $p$ which is computed in the $\operatorname{main}$ script by explicitly counting the parameters occurring in all MLPs that are used to parametrize the system architecture using PyTorch functions. The set of outcomes $O$ can then be read out from the logs which are managed using \textit{mlflow}~\cite{zaharia2018accelerating}.
{In order to estimate the compute in $\operatorname{flops}$, we use the general formula $c = p \cdot n_e \cdot n_t \cdot 1000$, which follows directly from the fact that each parameter contributes one floating-point operation per sample during the forward pass.}
\begin{algorithm}[t]\label{alg:pseudo}
\caption{Experimental setup using {input-state-output data}}\label{algo:submit}
\begin{algorithmic}[1]
\For{$A\!\in\!\{\text{unstructured~\eqref{eq:unstructured}},\text{input-affine~\eqref{eq:input_affine}},\text{pH~\eqref{eq:ph}}\}$}
    \For{$S\!\in\!\{\text{spring}~\eqref{eq:spring},\text{ball}~\eqref{eq:ball}, \text{motor}~\eqref{eq:motor}\}$}
        \For{$\operatorname{seed} \in \{0, \dots, N \} $}~~
            Set $\operatorname{seed}$
            \For{$n_e \in \{1, 2, 4, \dots, 2^{14} \}$}
                \For{$\tilde{d} \in \{2, 4, \ldots, 2^9\}$}
                    \For{$\tilde{n}_h \in \{ 2, 4, 8, 16 \}$}
                        \For{$\tilde{n}_d \in \{2, 4\}$}
                            \State Sample float $d \sim \{\tilde{d}/2, \dots, \tilde{T}\}$
                            \State {$n_t=\tfrac{1}{10}d$}
                            \State Sample integer $n_h \sim \{\tilde{n}_h/3, \dots, \tilde{n}_h\}$
                            \State Sample integer ${n}_d \sim \{\tilde{n}_d/2, \dots, \tilde{n}_d\}$
                            \State $\operatorname{nMAE}, \operatorname{nMSE}, p = \operatorname{main}(n_e, n_t, n_d, n_h, A, S)$
                            \State $c = p \cdot n_e \cdot n_t \cdot 1000$
                            \State Log $(\operatorname{nMAE}, \operatorname{nMSE}, c, p, d, A, S)$
                        \EndFor
                    \EndFor
                \EndFor
            \EndFor
        \EndFor
    \EndFor
\EndFor
\end{algorithmic}
\end{algorithm}

The total number of jobs submitted this way amounts to $63504$ which is executed on a high performance compute cluster with $268$ worker nodes, $17152$ Xeon Gold 6238R CPU cores and $931$ TB of parallel storage.

\subsection{{Model training using input-output data}}
{
We train the two model classes on two common benchmark tasks consisting of real data for nonlinear system identification:
\begin{itemize}
    \item The \textit{Wiener-Hammerstein system}~\cite{SCHOUKENS2017446} which is an electronic circuit with a diode-resistor nonlinearity;
    \item The \textit{Silverbox  system}~\cite{WigrenSchoukens13} which is an electronic implementation of a forced Duffing oscillator. 
\end{itemize}
A more detailed explanation of these systems and a link to the data repositories can be found in the appendix.
It is customary to use all possible lagged windows/output pairs in the entire dataset for one training iteration of NARX. Since mini-batch gradient descent is typically more computationally efficient, we use it to train NARX as well. Moreover, instead of a linear bypass in $f$, we used a residual connection \cite{he2016deep} of the form $f(x, u) = x + \tilde{f}(x, u)$ where $f$ is an MLP.}

{For training and validation, we vary the number of training data points by slicing subsets of the full benchmark datasets. If $n_\ell$ training data points are to be used and the full dataset consists of inputs $u_1, \dots, u_N$ and outputs $y_1, \dots, y_N$, then the subset considered is $u_1, \dots, u_{n_\ell}$ for the inputs and $y_1, \dots, y_{n_\ell}$ for the outputs. In addition, we adhere to the standard approach of separating out validation datasets, but we do not make use of the validation scores. Instead, the final score of each model on the test dataset is the score of the model at the last training iteration. For the Wiener-Hammerstein system, we consider $\tilde{n}_\ell \in\{ \lfloor 80000 \cdot 2^{-k}\rfloor: k=0, \ldots, 8\} =:L_{\mathrm{WH}}$ and $\tilde{n}_\ell \in\{\lfloor 64000 \cdot 2^{-k} \rfloor : k=0,\ldots, 8\} =: L_{\mathrm{Si}}$
for the Silverbox dataset. Then $n_\ell$ is chosen to be a random integer between $\frac{\tilde{n}_\ell}{2}$ and $\tilde{n}_\ell$ and the number of iterations for a job is $n_i:= m \cdot n_\ell$ where $m \in \{1, 2, \dots, 12\} =: M_S$ for SUBNET whereas $m \in \{1\}\cup\{ 2k : k=1,\ldots,15\} =: M_N$ for NARX. The reason for this difference in $m$ between NARX and SUBNET is that a training iteration for SUBNET is computationally more expensive and takes longer. Furthermore, the hidden dimensions $n_h$ considered are integers randomly sampled between $\frac{\tilde{n}_h}{2}$ and $\tilde{n}_h$ where $n_h \in \{ 2, 4, 8, 16, 32\} =: H_S$ for SUBNET and $n_h \in \{ 2, 4, 8, 16, 32, 64\} =: H_N$ for NARX. In addition, we consider the model depths $n_d \in  \{1, 2, 3\}:= D_{\mathrm{WH}}$ for Wiener-Hammerstein and $n_d \in \{1, 2\} =: D_{\mathrm{Si}}$ for 
Silverbox. }

{In the input-output case, we use $N = 6$ different random seeds. The expression $\operatorname{main}(n_i, n_\ell, n_d, n_h, A, S)$ in Line~28 of Algorithm~\ref{algo:submit2} is an abbreviation for submitting a job where a model with depth $n_d$ and hidden dimension $n_h$ is trained using $n_\ell$ data points for $n_i$ iterations on the dataset $D$ with architecture $A$.
This function then returns the metrics $\operatorname{nMAE}, \operatorname{nMSE}$ and the number of parameters of the model $p$ which is computed in the $\operatorname{main}$ script by explicitly counting the parameters occurring in all MLPs while excluding the parameters of the subspace encoder which is only called once per forward pass for SUBNET.}

{In order to compute the floating point operations for the NARX architecture, we use the formula $c = p \cdot n_i \cdot B$ where $p$ represents the number of parameters of the model. Since SUBNET forecasts for $T > 1$ steps during each training iteration, the formula $c = p \cdot n_i \cdot B \cdot T$ is used to compute the floating point operations used to train SUBNET. Finally, we use $T=130$ for the Wiener-Hammerstein system and $T=150$ for the Silverbox dataset. In each case the batch size $B=256$ is used. For NARX we use $n_u = n_y = 50$ and a past window of $50$ steps for the subspace encoder to compute the initial state for SUBNET.
}
\begin{algorithm}[t] \label{alg:pseudo2}
\caption{Experimental setup {using input-output data}}\label{algo:submit2}
\begin{algorithmic}[1]
\For{{$A\!\in\!\{\text{NARX},\text{SUBNET}\}$}}
    \For{{$S\!\in\!\{\text{Wiener-Hammerstein},\text{Silverbox}\}$}}
        \State {$D_{\mathrm{WH}} := \{ \lfloor 8\cdot10^4\cdot 2^{-k}\rfloor: k=0, \ldots, 8\}$
        \State $D_S = \{ 64\cdot 10^3 \cdot 2^{-k}  : k=0,\ldots, 8\}$}
        \State {$H_S := \{2, 4, 8, 16, 32\}; H_N = \{2, 4, 8, 16, 32, 64\}$}
        \State {$M_S := \{1, 2, \dots, 12\}$
        \State $M_N:=\{1\}\cup\{ 2k : k=1,\ldots,15\}$}
        \State {$H_S := \{2, 4, 8, 16 32\}$; $H_N := \{2, 4, 8, 16, 32 64\}$}
        \State {$D_{\mathrm{WH}} := \{1, 2, 3\}$; $D_S := \{1,2\}$}
        \If{{A == SUBNET}}
        \State {$M := M_S$; $H :=H_S$}
        \Else
        \State {$M:= M_N$; $H:= H_N$}
        \EndIf
        \If{{D == Wiener-Hammerstein}}
        \State {$D := M_{\mathrm{WH}}$; $L :=L_{\mathrm{WH}}$}
        \Else
        \State {$D := M_{\mathrm{Si}}$; $L :=L_{\mathrm{Si}}$}
        \EndIf
        \For{{$\operatorname{seed} \in \{0, \dots, N \} $}~~}
            {Set $\operatorname{seed}$}
            \For{{$m \in M$}}
                \For{{$\tilde{n}_\ell \in L$}}
                    \For{{$\tilde{n}_h \in H$}}
                        \For{{$n_d \in D$}}
                            \State {Sample integer $n_\ell \sim \{ \tilde{n}_\ell/2, \dots, \tilde{n}_\ell \} $}
                            \State {$n_i := m \cdot n_\ell$}
                            \State {Sample integer $n_h \sim \{\tilde{n}_h/2, \dots, \tilde{n}_h\}$}
                            \State {$\operatorname{nMAE}, \operatorname{nMSE}, p = \operatorname{main}(n_i, n_\ell, n_d, n_h, A, S)$}
                            \State {$c = p \cdot n_i \cdot B$}
                            \If{{A == SUBNET}}
                            \If {{S == Wiener-Hammerstein}}
                            \State {$c = c \cdot 130$}
                            \Else 
                            \State {$c = c \cdot 150$}
                            \EndIf
                            \EndIf
                            \State {Log $(\operatorname{nMAE}, \operatorname{nMSE}, c, p, n_\ell, A, S)$}
                        \EndFor
                    \EndFor
                \EndFor
            \EndFor
        \EndFor
    \EndFor
\EndFor
\end{algorithmic}
\end{algorithm}

\subsection{Approximation of Neural Scaling Laws}
The NSLs are fitted semi-automatically. We first manually construct a rough approximation of the NSLs by specifying a piecewise affine function at $\log-\log$ scale and then automatically fit the NSLs starting with this approximation. We find that this initial guess is necessary since the problem of fitting NSLs is non-convex.
To fit the NSLs, we use gradient optimization in PyTorch~\cite{Pytorch}. More precisely, we parameterize a module with parameters representing the values $\log \alpha, \log \beta, \delta, \log \sigma_i, \varphi_i, \delta_i$ for $i = 1, \dots, b$ in the representation~\eqref{eq:broken_NSL}. The value $\varphi_i$ is capped to be higher than $0.2$ by replacing it with $\max(\varphi_i, 0.2)$   whenever the module is called. This ensures numerical stability during fitting.
We point out that learning $\log \beta, \log \sigma_i$ instead of $\beta, \sigma_i$ is necessary since the latter are expected to have a high order of magnitude. It would otherwise {take an impractically long time} to successfully complete the fitting, because gradient steps are typically very small.

Assume that the {lower envelope $E_O$ given by \eqref{envelope}}, the interpolation points $\tilde{r}_1, \dots, \tilde{r}_K$, and the evaluation of $E_O$ on those points $\tilde{e}_1, \dots, \tilde{e}_K$ are given as in Definition~\ref{def: broken nsl} and we choose the number of interpolation points $K = 100$ for the fitting of all sets of outcomes. 
The \textit{ADAM} optimizer \cite{kingma2017adammethodstochasticoptimization} is used to train the module for a few {hundred} iterations, in which the whole dataset is passed at each iteration.
{The loss function being minimized is \[
\mathcal{L} = \frac{1}{K}\sum\limits_{k=1}^K (\beta \cdot \max(\hat{e}_k-\tilde{e}_k, 0) + \max(\tilde{e_k}-\hat{e}_k, 0))
\]
where $\tilde{e}_1, \dots, \tilde{e}_K$ are given as in Definition~\ref{def: broken nsl} with $\hat{e}_1, \dots, \hat{e}_K$ being the current approximations. We use $\beta = 5$ to give a higher weight to errors where the approximated NSL lies above the lower envelope relative to errors where it lies below the lower envelope.}

\section{Verification of neural scaling laws for input-state-output data}\label{sec:verification_ios}

In this section, we show that several types of {broken} NSLs can be observed in the context of learning-based system identification {using input-state-output data}. We show that different identification methods yield reproducible NSLs across a {different} system types and prior physical knowledge levels. 
An overview of all the experiments and approximated NSLs is presented in the Appendix. In the following subsections, we take a closer look at a few figures that represent the general trends we noticed by inspecting all the results of the different experiments.

\subsection{Compute neural scaling laws}\label{sec:compute}
In this section, we consider the compute NSL by considering the total amount of compute $c$ in terms of $\operatorname{flops}$ used to train the model.
Since the other parameters (data size and model size) influence the compute NSL, these parameters need to be varied simultaneously. 
\begin{figure}[!ht]
    \centering
    \includegraphics[width=1\columnwidth]{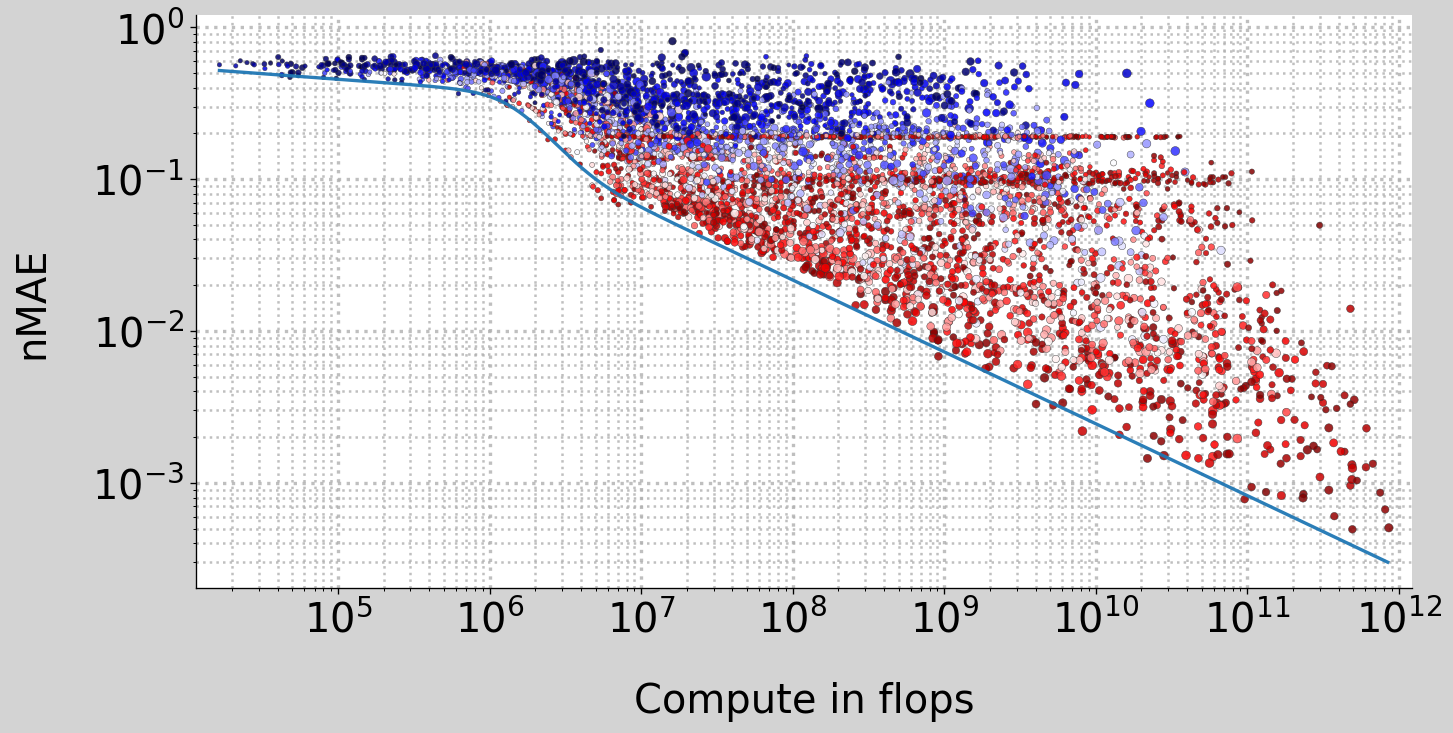}\\[0.5em]
    \includegraphics[width=.65\columnwidth]{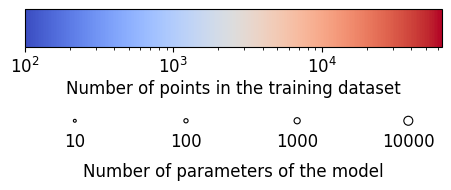}
    \caption{Compute-$\operatorname{nMAE}$ 
    NSL for the input-affine architecture on data generated from the ball system { which has a margin $M=0.04$ and $b=2$ breaks.}}
    \label{fig:compute-nmae-aff-ball}
\end{figure}
Figure~\ref{fig:compute-nmae-aff-ball} shows the compute-$\operatorname{nMAE}$ NSL for the input-affine architecture on the ball system using one break.
The $x$-axis represents the compute in $\operatorname{flops}$ used to train the model (resource being scaled) and the $y$-axis represents the error (measured as $\operatorname{nMAE}$).
Individual experiments are represented by points in the plot.
{The point color indicates the size of the dataset (in terms of individual data points)}  whereas the point size represents the model size (in terms of total parameters). The approximate NSL given by \eqref{eq:broken_NSL} determined as described in the beginning of Section~\ref{sec:approximation} is represented by the blue curve.
 It is given by the formula 
{\begin{align*}
L(c) = 0 + 1.1c^{-0.073}&\left[1 + \left(\frac{c}{1.4\times 10^{6}}\right)^{ \frac{1}{0.2} }\right]^{-1.1 \cdot 0.2}\left[1 + \left(\frac{c}{4.6\times 10^{6}}\right)^{ \frac{1}{0.2} }\right]^{0.68 \cdot 0.2}
.\end{align*}}
We can see that the performance of models does not increase with more compute until the plateau of random guessing is being overcome at roughly $10^6 \ \operatorname{flops}$. The curve then follows a small segment with steep slope. It then enters a large segment where it decays linearly on $\log$-$\log$ scale (which means that it approximately satisfies a power law) until approximately $10^{12} \operatorname{flops}$ which is the highest compute budget we considered. Hence, the NSL shown in Figure~\ref{fig:compute-nmae-aff-ball} has $b=2$ breaks.

We note that within the range of $10^7 \ \operatorname{flops}$ up to approximately $10^{12} \ \operatorname{flops}$, there is highly predictable scaling with the given amount of data and model parameters that we consider. We can see that the smallest points in the plot lie far away from the envelope, which shows that models need to have a minimum amount of parameters to yield compute-optimal performance. We can also see that it is necessary to have a minimum amount of data available to obtain scaling, which is indicated by points with low errors typically having a stronger red tone. Moreover, it is interesting to note that the point corresponding to the highest compute being used lies close to the lower envelope $E_O$. This demonstrates that training models for as much as the unexpectedly high number $n_e = 10^{14}$ of total epochs {can still decrease the evaluation error.} 
By inspecting the overview in Figure~\ref{fig:comp-nmae}, we can see that the compute-$\operatorname{nMAE}$  NSLs look very similar for all architectures and datasets we consider. This indicates that the compute necessary to obtain a certain $\operatorname{nMAE}$-value is highly predictable.

\begin{figure}[!ht]
    \centering
    \includegraphics[width=1\columnwidth]{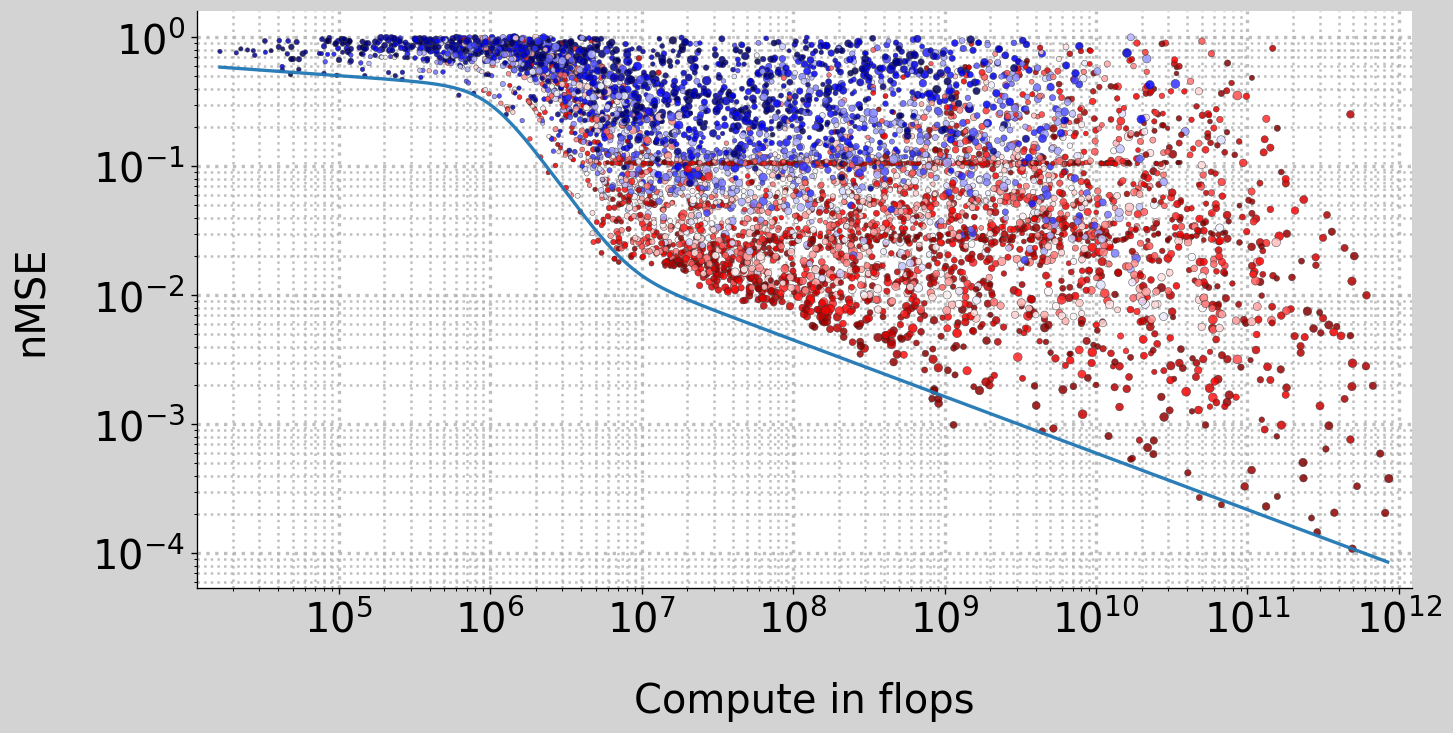}\\[0.5em]
    \caption{Compute-$\operatorname{nMSE}$
    NSL for the input-affine architecture on data generated from the ball system {which has a margin $M=0.09$ and $b=2$ breaks}.}
    \label{fig:compute-nmse-aff-ball}
\end{figure}
Figure~\ref{fig:compute-nmse-aff-ball} shows the corresponding  compute-$\operatorname{nMSE}$ NSL for the input-affine architecture on the {ball} system. 
The $x$-axis again represents the compute in $\operatorname{flops}$ used to train the model and the $y$-axis represents the error (measured as $\operatorname{nMSE}$ or $\operatorname{nMAE}$). The NSL in Figure~\ref{fig:compute-nmse-aff-ball} is given by the formula
{\begin{align*}
L(c) = 0 + 1.3c^{-0.084}&\left[1 + \left(\frac{c}{10^{6}}\right)^{ \frac{1}{0.2} }\right]^{-1.6 \cdot 0.2}\left[1 +\left(\frac{c}{8.1\times 10^{6}}\right)^{ \frac{1}{0.2} }\right]^{1.2 \cdot 0.2}.
\end{align*}}
We observe $b=2$ breaks at similar compute budgets for both error measures. We note, however, that at more than $10^{10} \operatorname{flops}$ there are fewer points close to the envelope for $\operatorname{nMSE}$. Inspecting the overview in Figures~\ref{fig:comp-nmae}~and~\ref{fig:comp-nmse}, one observes, in fact, that the $\operatorname{nMAE}$ NSLs tend to have fewer outliers and that there are more points located at the envelope. This suggests that $\operatorname{nMAE}$ could be a more suitable metric than $\operatorname{nMSE}$ in the context of compute NSLs.

We conclude that compute NSLs allow to predict the amount of compute (in $\operatorname{flops}$) needed to obtain a certain accuracy within a very high range of possible configurations.

\subsection{Data neural scaling laws}\label{sec:data}

\begin{figure}[!ht] 
    \centering
    \includegraphics[width=1\columnwidth]{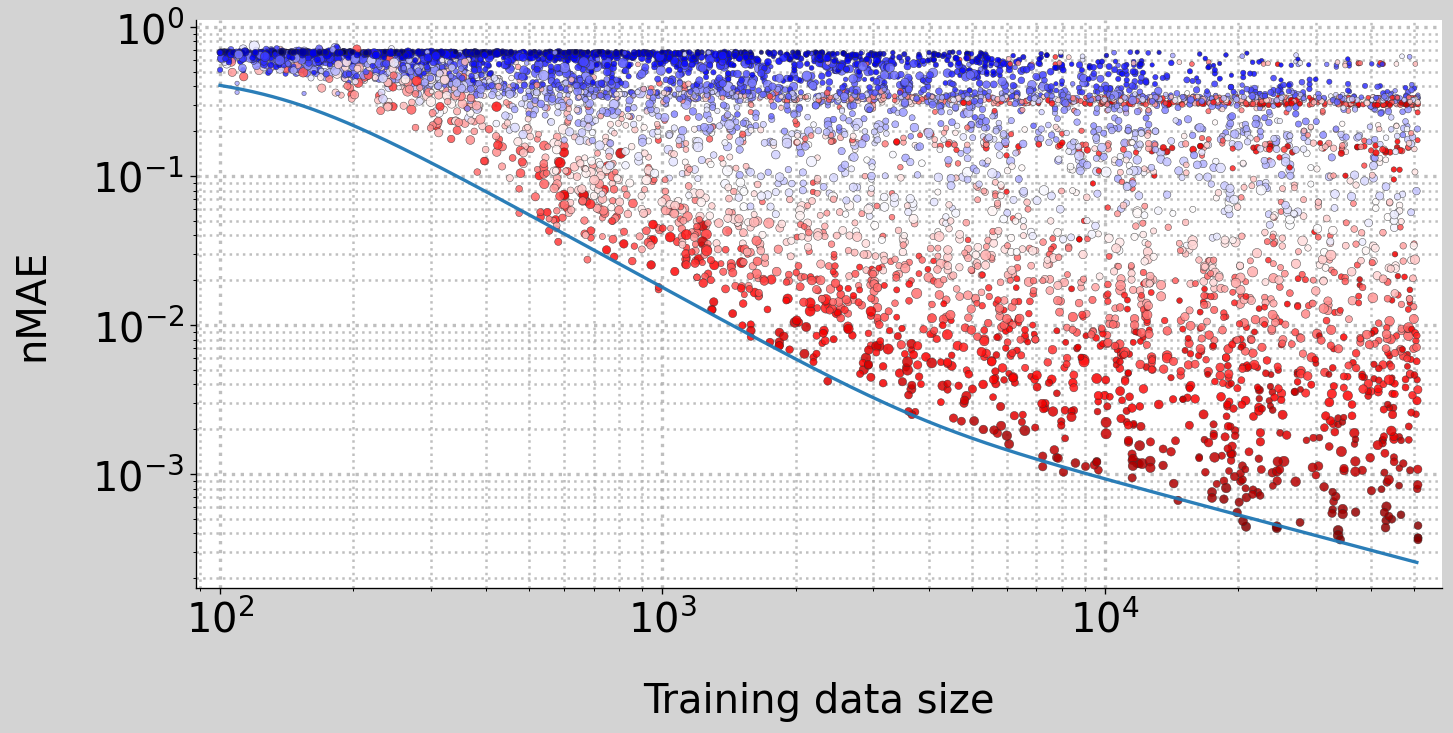}\\[0.5em]
    \includegraphics[width=0.65\columnwidth]{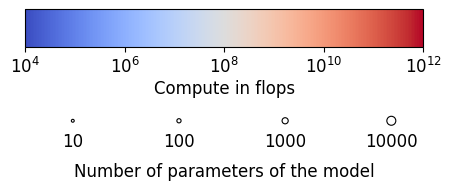}
    \caption{Data-$\operatorname{nMAE}$  NSL for the unstructured architecture on the motor system {which has a margin $M=0.07$ and $b=2$ breaks}.}
    \label{fig:data-nmae-uns-motor}
\end{figure}

\begin{figure}[!ht] 
    \centering
    \includegraphics[width=1\columnwidth]{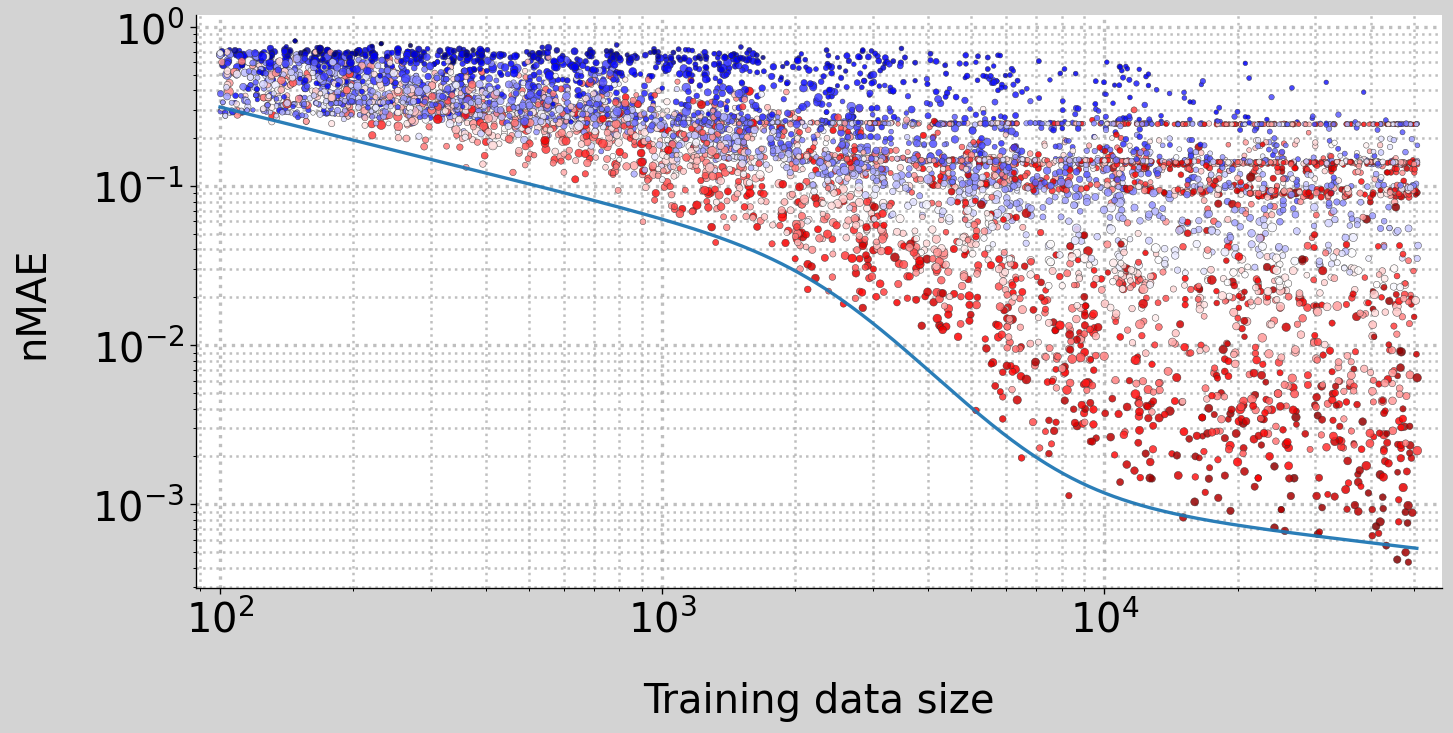}\\[0.5em]
    \caption{Data-$\operatorname{nMAE}$  NSL for the input-affine architecture on the spring system {which has a margin $M=0.07$ and $b=2$ breaks.}}
    \label{fig:data-nmae-aff-spring}
\end{figure}

In this section, we consider data NSLs, which are obtained by treating the {number of points} in the dataset as the resource being varied. 
Each point in the Figures~\ref{fig:data-nmae-aff-spring} and \ref{fig:data-nmae-uns-motor} represents an element of the set of outcomes $O$ for different training configurations. Here, the color of the points represents the amount of compute measured in $\operatorname{flops}$ and the number of model parameters is represented by the point size. 
The $x$-axis represents the total {number of points $n_\ell$ in the training dataset. In the case of input-state-output data $n_\ell = n_t \cdot 1000$ while it is logged as an input to the main function for input-output data.}

The data-$\operatorname{nMAE}$ NSL {for unstructured architecture on the motor system} is shown in Figure~\ref{fig:data-nmae-uns-motor}{,} it has $b=2$ breaks and is given by the formula
{\begin{align*}
   L(d) = 0 + 0.67C^{-0.079}&\left[1 + \left(\frac{d}{1.4\times 10^{2}}\right)^{ \frac{1}{0.2} }\right]^{-1.5 \cdot 0.2}\left[1 + \left(\frac{d}{3.9\times 10^{3}}\right)^{ \frac{1}{0.2} }\right]^{0.84 \cdot 0.2}.
\end{align*}  }
We notice again that a minimum amount of data is necessary to overcome the plateau of random guessing which is followed by a small segment of fast progression and then ends with a long segment of smooth scaling. Similar to the compute NSL, it can be observed that the remaining two resources compute and model size play a significant role. More precisely, we observe that points close to the envelope tend to be red, especially for the higher amounts, which indicates that enough compute needs to be used to follow the data NSL. Moreover, we observe that very small points are typically not close to the envelope, which indicates that models that are too small cannot optimize the performance for a given amount of data. In addition, we observe a smooth scaling approximately from {$n_\ell = 8000$} to ${n_\ell = 64000}$, which indicates that the amount of data necessary to observe the desired performance can be predicted considering experiments with smaller amounts of data in this regime. Since the NSL perfectly mimics the envelope, these predictions can be considered to have a high confidence.

Similar behavior can be noticed in Figure~\ref{fig:data-nmae-aff-spring} that shows the data-$\operatorname{nMAE}$ NSL for the input-affine architecture on the spring dataset which is given by the formula
{\begin{align*}
    L(d) = 0 + 7.5C^{-0.69}&\left[1 + \left(\frac{d}{2.4\times 10^{3}}\right)^{ \frac{1}{0.2} }\right]^{-2.4 \cdot 0.2}\left[1 + \left(\frac{d}{7.6\times 10^{3}}\right)^{ \frac{1}{0.2} }\right]^{2.8 \cdot 0.2}.
\end{align*}}However, compared to Figure~\ref{fig:data-nmae-uns-motor}, we can see that there {is a break approximately at the relatively high dataset size $n_\ell = 10^4$.} {This means that it is not possible to predict the performance achievable with more than $10^4$ data points solely from experiments with less than $10^4$ data points.}
{This} contrasts with the observations made for the compute NSLs in Section~\ref{sec:compute} and the observations made for the unstructured architecture applied to the motor system shown in Figure~\ref{fig:data-nmae-uns-motor}.
By inspecting Figures~\ref{fig:data-nmae}~and~\ref{fig:data-nmse}, we can moreover see that data NSLs look very similar for the different identification architectures considered given a fixed dataset which means that they generalize well. 

We conclude that data NSLs can be used to predict the error of a model at a given amount of data. However, these data NSLs can suffer from the issue that {breaks can occur at relatively large dataset sizes. This means that while there is a good chance that the extrapolation of the NSL yields an accurate estimate of the data needed to achieve a given desired performance, there is no guarantee.} This phenomenon can also be seen in Figures~\ref{fig:data-nmae}~and~\ref{fig:data-nmse}, for example for the $\mathrm{nMSE}$-metric on the ball system.

\subsection{Model neural scaling laws}\label{sec:model}
In this section, we consider model NSL which are obtained by treating the number of parameters $p$ of the model as the resource being varied. 
\begin{figure}[!ht] 
    \centering
    \includegraphics[width=1\columnwidth]{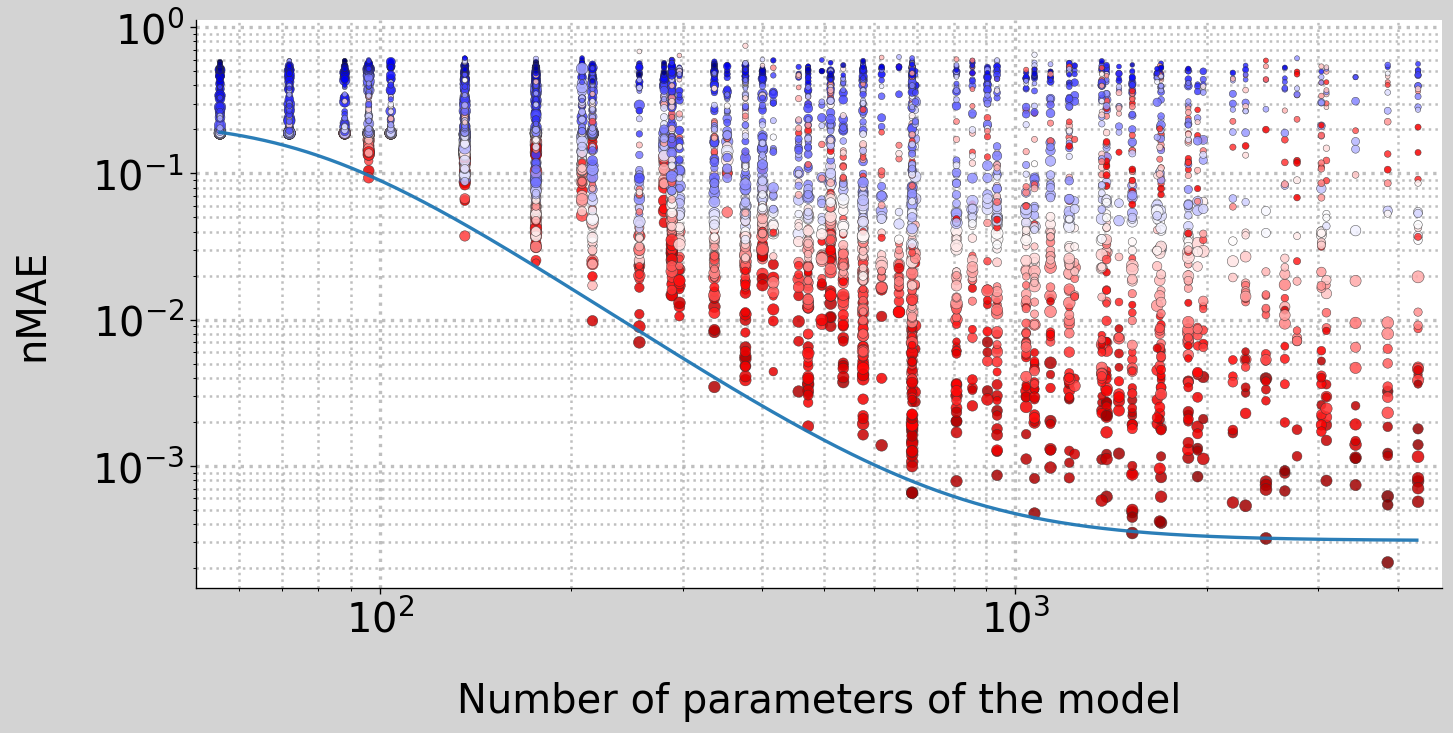}\\[0.5em]
    \includegraphics[width=0.65\columnwidth]{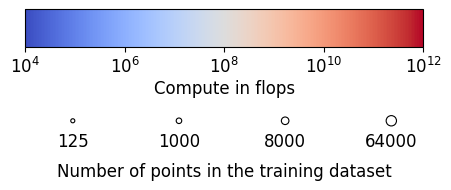}
    \caption{Model-$\operatorname{nMAE}$ NSL for the port-Hamiltonian architecture on the ball system {which has a margin $M=0.11$ and $b=1$ breaks}.}
    \label{fig:model-nmae-ph-ball}
\end{figure}

Figure~\ref{fig:model-nmae-ph-ball} shows the model NSL for the pH~architecture for the ball system. Here, the $x$-axis represents the number of parameters of the model. The point size indicates the dataset size and the color represents the total compute used to train the model.
This NSL can be well represented using an approximation with a single break and is given by the formula
{\begin{align*}
    L(p) = 0.00031 + 0.12p^{0.19}\left[1 + \left(\frac{p}{76}\right)^{ \frac{1}{0.2} }\right]^{-3.1 \cdot 0.2}.
\end{align*}}
We can see that scaling the model further than approximately $10^3$ parameters does not yield better results. Since this is also the case for our configuration, we infer that scaling the model parameter size is not the major bottleneck {because a model with $10^3$ parameters can be considered small when executed on standard modern consumer hardware}. However, it is important to note that models which are too small do not achieve good performance. We can also see that runs with insufficient compute or insufficient data do not achieve good performance at a given model size. This is not surprising since we have also observed in Sections~\ref{sec:compute}~and~\ref{sec:data} that all three resources are important to consider if one wants to improve the performance for a given resource.

We conclude that the model size needed to obtain a given error can be determined using NSL. We also noticed that scaling up the model was not a computational challenge in our experiments. However, it is important to scale the number of parameters of the model together with data and compute to obtain optimal results, see Figure~\ref{fig:comp-nmae}.

\section{Verification of neural scaling laws for input-output data}\label{sec:verification_io}
{In this section we investigate NSLs for input-output data. We will focus on the nMAE-NSLs, the corresponding nMSE-NSLs can be found in the appendix.}

\subsection{{Compute neural scaling laws}}

\begin{figure}[!ht] 
    \centering
    \includegraphics[width=1\columnwidth]{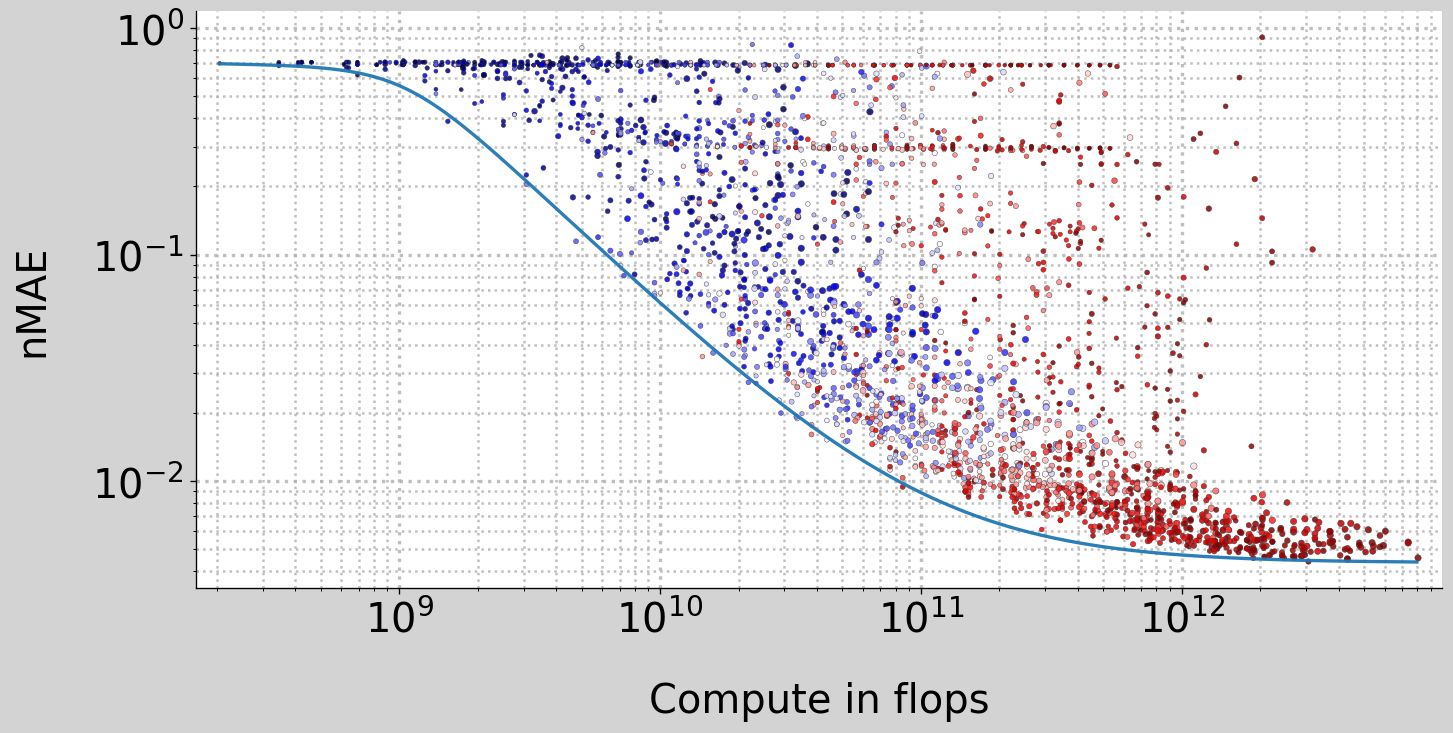}\\[0.5em]
    \caption{{Compute-$\operatorname{nMAE}$ NSL for the SUBNET architecture on the Silverbox system {which has a margin $M=0.02$ and $b=1$ break}.}}
    \label{fig:compute-subnet-sb}
\end{figure}
{
In Figure~\ref{fig:compute-subnet-sb}, we can see the compute NSL of the SUBNET architecture on the Silverbox system.
This NSL is represented by the formula \[
L(c) = 0.0054 + 0.13c^{-0.0029}\left[1 + \left(\frac{c}{4\times 10^{8}}\right)^{ \frac{1}{0.2} }\right]^{-0.73 \cdot 0.2}
\]
and has one break.
We can clearly see that this NSL has a very similar structure as the compute NSLs for input-state-output data system identification. However, one notable difference to most of the latter ones is that it shows strong saturation for more than $10^{12}\, \mathrm{flops}$. This indicates that using more compute to train SUBNET on this dataset will most likely not yield better results. That is not a surprise, since the dataset size is bounded by the amount of measurements taken in reality, whereas the dataset sizes for the simulation datasets considered in the input-state-output case were chosen so that no strong saturation happens.
}

\begin{figure}[!ht] 
    \centering
    \includegraphics[width=1\columnwidth]{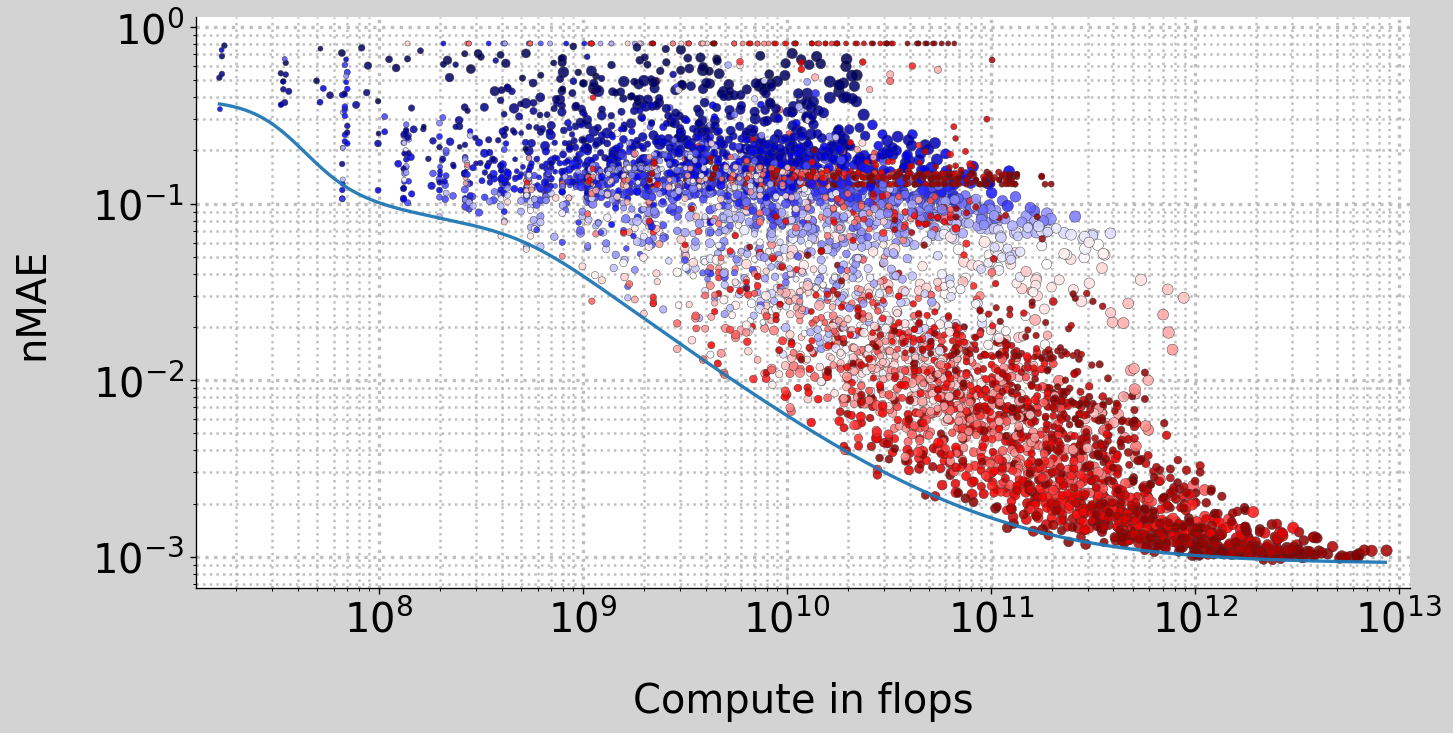}\\[0.5em]
    \caption{{Compute-$\operatorname{nMAE}$ NSL for the NARX architecture on the Wiener-Hammerstein system which has a margin $M=0.05$ and $b=3$ breaks.}}
    \label{fig:compute-narx-wh}
\end{figure}
{Figure~\ref{fig:compute-narx-wh} represents the compute NSL of the NARX architecture on the Wiener-Hammerstein system.
This NSL is represented by the formula \begin{align*}
    L(c) = 0.00091 + 0.34c^{0.0081}&\left[1 + \left(\frac{c}{3.7\times 10^{7}}\right)^{ \frac{1}{0.2} }\right]^{-4.8 \cdot 0.2}\\
    &\left[1 + \left(\frac{c}{4.8\times 10^{7}}\right)^{ \frac{1}{0.2} }\right]^{4.6 \cdot 0.2}\\
    &\left[1 + \left(\frac{c}{5.3\times 10^{8}}\right)^{ \frac{1}{0.2} }\right]^{-0.66 \cdot 0.2}
\end{align*}
and has three break which all occur in the region from $10^7\, \mathrm{flops}$ to $10^9\, \mathrm{flops}$.
We do again observe a strong saturation at more than $10^{12}\, \mathrm{flops}$, indicating that it is unlikely to improve the performance by using more resources to train the models.
}

\subsection{{Data neural scaling laws}}
\begin{figure}[!ht] 
    \centering
    \includegraphics[width=1\columnwidth]{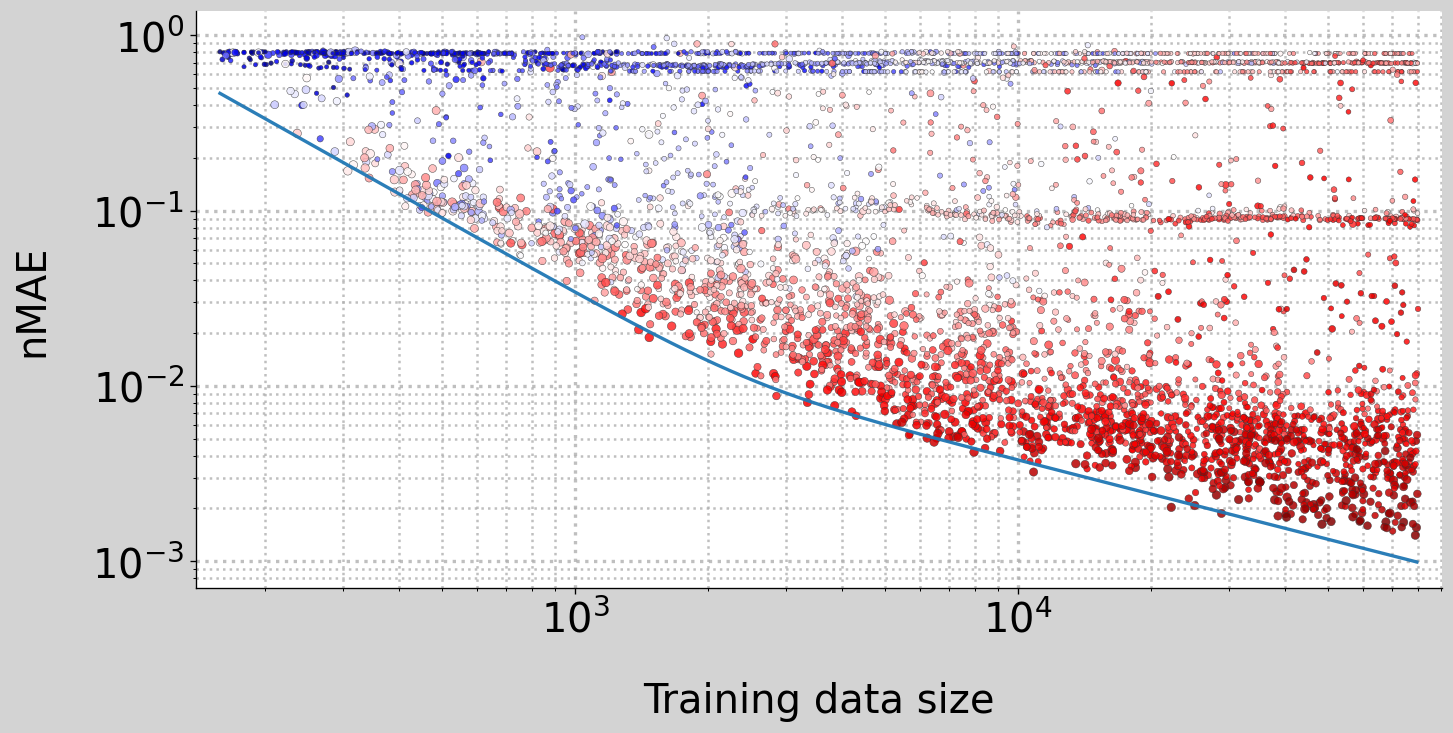}\\[0.5em]
    \caption{{Data-$\operatorname{nMAE}$ NSL for the SUBNET architecture on the Wiener-Hammerstein  system which has a margin $M=0.05$ and $b=1$ break.}}
    \label{fig:data-subnet-wh}
\end{figure}
{
In Figure~\ref{fig:data-subnet-wh}, we can see the data-$\operatorname{nMAE}$ NSL for SUBNET on the Wiener-Hammerstein system.
The functional form of the NSL is given by the formula \begin{align*}
    L(d) = 0 + 6.2\times 10^{2}d^{-1.4}\left[1 + \left(\frac{d}{2.5\times 10^{3}}\right)^{ \frac{1}{0.2} }\right]^{0.77 \cdot 0.2}
\end{align*}
and has one break.
Since we have seen strong saturation of the compute-nMAE NSL on the Wiener-Hammerstein system for the NARX architecture (we do in fact also for SUBNET as one can see Figure~\ref{fig:comp-io-nmae}), we suspect that the data-NSL is not supposed to show saturation. We can in fact see that this is the case, because the NSL shown in Figure~\ref{fig:data-subnet-wh} exhibits smooth scaling at the highest data sizes we considered (size of the full Wiener-Hammerstein system itself). This indicates that it might be relatively easy to identify the system with more accuracy by enlarging the dataset.
\begin{figure}[!ht] 
    \centering
    \includegraphics[width=1\columnwidth]{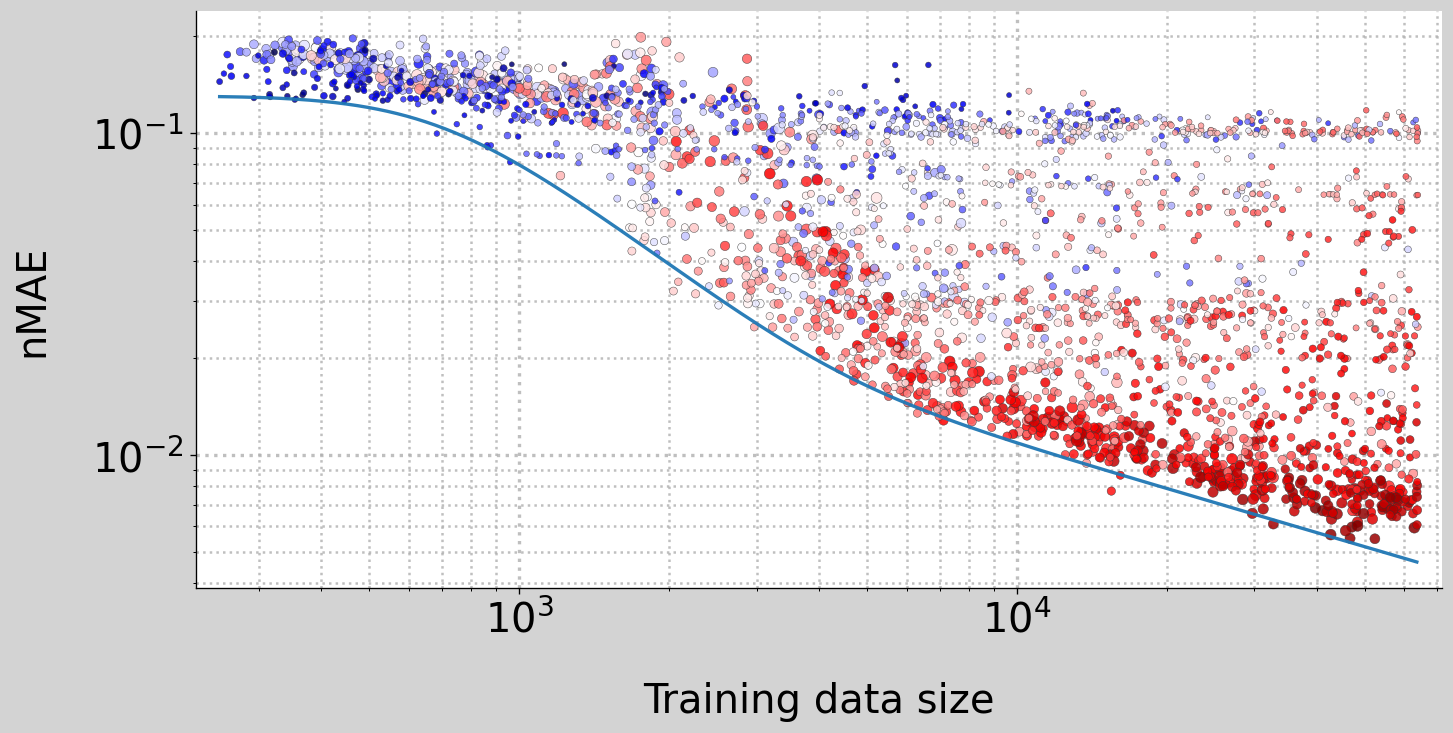}\\[0.5em]
    \caption{{Data-$\operatorname{nMAE}$ NSL for the NARX architecture on the Silverbox system which has a margin $M=0.05$ and $b=2$ breaks.}}
    \label{fig:data-narx-sb}
\end{figure}
}
{
In Figure~\ref{fig:data-narx-sb}, we can see the data-$\operatorname{nMAE}$ NSL for the NARX architecture on the Silverbox system. It is represented by the formula
\begin{align*}
    L(d) = 0 + 0.12d^{0.019}&\left[1 + \left(\frac{d}{6.9\times 10^{2}}\right)^{ \frac{1}{0.2} }\right]^{-1.2 \cdot 0.2}\left[1 + \left(\frac{d}{4.4\times 10^{3}}\right)^{ \frac{1}{0.2} }\right]^{0.71 \cdot 0.2}
\end{align*}
and has two breaks. We do again observe that the smooth scaling regime continues towards the size of the full Silverbox system, indicating that enlarging the dataset might greatly improve the performance of the system identification task. 
Moreover, we can use the slopes observed at $10^4$ in Figures~\ref{fig:data-narx-sb}~and~\ref{fig:data-subnet-wh} to anticipate how large the dataset has to be in order to achieve a given performance margin for the system identification task.
}
\subsection{{Model neural scaling laws}}
\begin{figure}[!ht] 
    \centering
    \includegraphics[width=1\columnwidth]{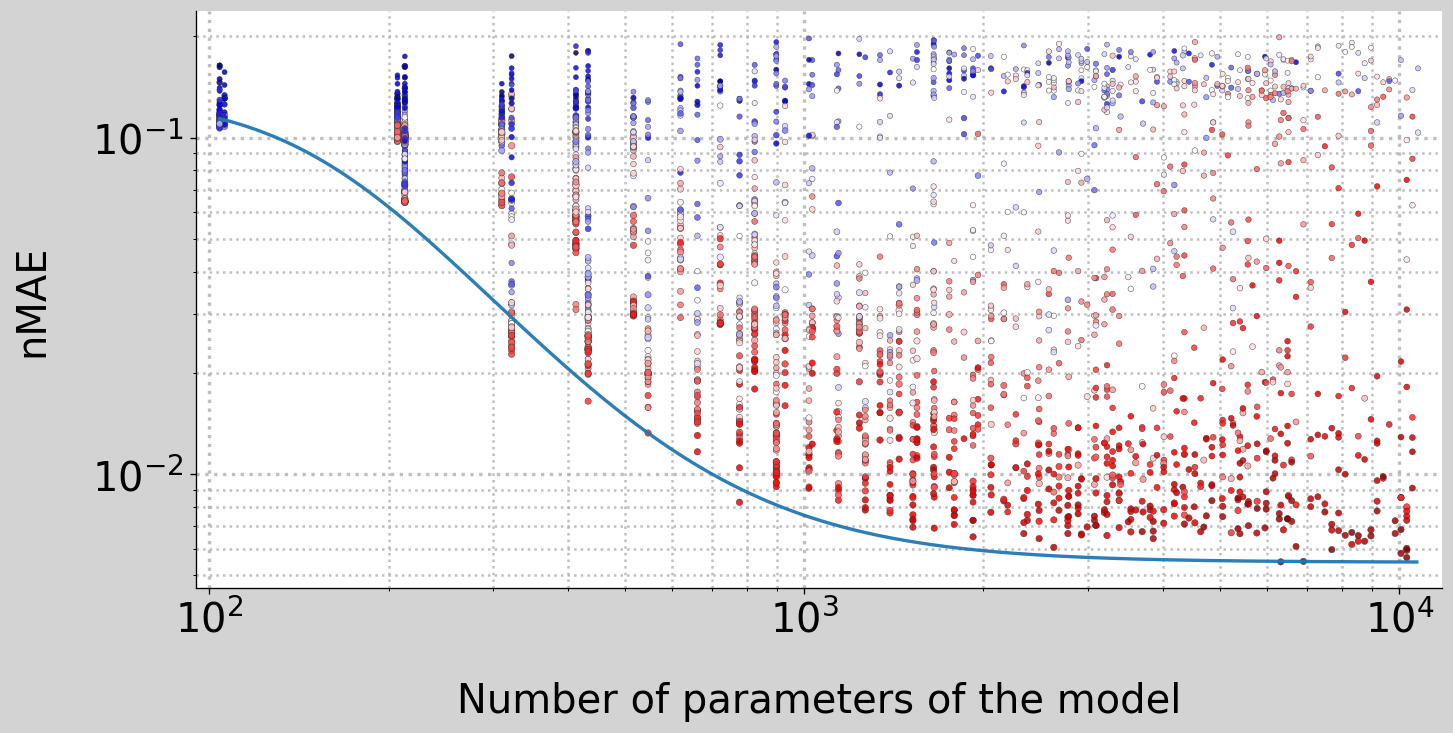}\\[0.5em]
    \caption{{Data-$\operatorname{nMAE}$ NSL for the NARX architecture on the Silverbox system which has a margin $M=0.09$ and $b=1$ break.}}
    \label{fig:model-narx-sb}
\end{figure}
{In Figure~\ref{fig:model-narx-sb}, we can see the model NSL for the NARX architecture on the Silverbox system. This NSL looks very similar to the corresponding model-NSLs for input-state-output data system identification. The NSL is represented by the function
\[
L(p) = 0.0055 + 0.14p^{-0.014}\left[1 + \left(\frac{p}{1.5\times 10^{2}}\right)^{ \frac{1}{0.2} }\right]^{-2.2 \cdot 0.2}.
\]
}
{
The NSL saturates at the relatively small amount of $10^3$ parameters. This shows that for the considered dataset sizes, we would not have gotten better results using larger models. We can moreover also see in Figure~\ref{fig:comp-io-nmae} that the compute-nMAE NSL for the NARX architecture on the Silverbox dataset is also highly saturated.}

{
On the contrary, we can see in Figure~\ref{fig:data-narx-sb} that adding more data is expected to yield better results, because the data NSL does not saturate. Since the compute and the model NSLs are highly saturated for this task and the data NSL not, we can conclude that adding more data is significantly more likely to yield better performance gains than scaling up the other two resources compute and model size.
}
\section{Conclusion}
Neural Scaling Laws (NSLs) quantify the dependence between accuracy, dataset size, model complexity, and compute available. Here, we show that NSLs can also be verified in the context of learning-based system identification for various examples and system architectures. Early trends in the scaling curve allow us to anticipate the data, compute and model complexity required for higher accuracy.

Furthermore, we observed that in system identification, the amount of data is more important than model size, which is in contrast to the results in the natural language processing community. This means that a larger model does not necessarily yield better results, as indicated by the model NSLs that are often saturating with relatively few parameters. Moreover, system models must be trained for a very large number of iterations, since training for up to $2^{14}$ epochs still produces improvements in the performance  of the identified system model for different datasets, dataset sizes, system architectures, and model sizes.

\bibliographystyle{IEEEtran}
\bibliography{refs}

\begin{appendices}

\section{Collection of example systems}\label{app:ex}

We use the same examples as in \cite{cherifi2025}.

\subsection{Spring system}
We consider the interconnection of two nonlinearly damped mass-spring systems from \cite{Neary2023} with 
\begin{equation}
     H(x) = \sum_{i=1}^2\frac{k_i q_i^2}{2} + \frac{p^2_i}{2m_i}, \: R(x) = \begin{bmatrix}
        0 & 0 & 0 & 0 \\
        0 & \frac{b_1p_1^2}{m_1^2} & 0 & 0 \\
        0 & 0 & 0 & 0 \\
        0 & 0 & 0 & \frac{b_2p_2^2}{m_2^2}
    \end{bmatrix},
    \label{eq:spring}
\end{equation}
\begin{equation*}
 B(x) = \begin{bmatrix}
        0 & 1&0&0\\
        0& 0 &0&1
    \end{bmatrix}^\top, \:  J(x) = \begin{bmatrix}
        0 & 1 & 0 & 0 \\
        -1 & 0 & 1 & 0 \\
        0 & -1 & 0 & 1 \\
        0 & 0 & -1 & 0
    \end{bmatrix},
\end{equation*}
with $x = \begin{bmatrix}
        q_1 &
        p_1 &
        q_2 &
        p_2
    \end{bmatrix}^\top$.
\begin{table}[h!]
\centering
\renewcommand{\arraystretch}{1.5}
\begin{tabular}{|c |c |c |c |}
\toprule
\hline
symbol& meaning & value & unit \\
\hline\hline
$m_1$ & mass 1 & {0.5} & $[{\rm kg} ]$ \\ $m_2$ & mass 2 & {0.75} & $[{\rm kg} ]$ \\ 
$b_1,b_2$ & damping coefficient & {0.5} & - \\
$k_1$ & spring constant 1 & {2} & $[{\rm N/m}]$ \\
$k_2$ & spring constant 2 & {0.2} & $[ {\rm N/m}]$ \\
\hline
\bottomrule
\end{tabular}
\caption{Mass spring system parameters}
\label{tab:mass_spring}
\end{table}
    The parameter values are shown in Table~\ref{tab:mass_spring}.
We will refer to this system as the \textit{spring system}.

\subsection{Ball system}
The second system we consider is magnetically levitated iron ball system from \cite{beckers2022gaussian}, which is referred to as the \textit{ball system}. In that case
\begin{align}
\label{eq:ball}
    H(x) &= \frac{{x_2}^2}{2m} + \frac{1}{2} \frac{x_3^2}{L(x_1)},& B(x) &= \begin{bmatrix}
    0 &
    0 &
    1
\end{bmatrix}^\top ,\\
    J(x) &= \begin{bmatrix}
        0 & 1 & 0 \\
        -1 & 0 & 0 \\
        0 & 0 & 0
    \end{bmatrix},& R(x) &= \begin{bmatrix}
        0 & 0 & 0 \\
        0 & c|x_2| & 0 \\
        0 & 0 & \frac{1}{R}
    \end{bmatrix},\nonumber
\end{align}
where $x = \begin{bmatrix}
    x_1 &
    x_2 &
    x_3
\end{bmatrix}^\top$ with $x_1$ representing the vertical position of the ball, $x_2$ its momentum and $x_3$ the magnetic flux.
\begin{table}[!h]
\centering
\renewcommand{\arraystretch}{1.5}
\begin{tabular}{|c |c |c |c |}
\toprule
\hline
symbol& meaning & value & unit \\
\hline\hline
$m$ & ball mass & 0.1 & $[{\rm kg}]$ \\  
$L(x_1)$ & inductivity at height $x_1$ & $(0.1+x_1^2)^{-1}$ & $[{\rm H}]$ 
\\
$R$ & electrical resistance  & 0.1 & $[\Omega]$ \\
$c$ & drag coefficient  & 1 & -\\
\hline
\bottomrule
\end{tabular}
\caption{Ball system parameters}
\label{tab:ironball}
\end{table}
The parameter values are shown in Table~\ref{tab:ironball}. 

\subsection{Motor system}
\label{sec:motor}
The third system we consider is the permanent magnet synchronous motor, see \cite{vu2023port, spirito2024structure}, which is referred to as the \textit{ motor system}.

The system structure is 
\begin{align}
\label{eq:motor}
    H(x) &= \frac{\varphi_d^2}{2L} +  \frac{\varphi_q^2}{2L} + \frac{p^2}{2J_m}, \quad  B = \begin{bmatrix}
    1 & 0 &0 \\
    0 & 1 &0
\end{bmatrix}^\top\\
    J(x) &\!=\! \begin{bmatrix}
        0 & 0 & \varphi_q \\
        0 & 0 & - \varphi_d - \Phi \\
        -\varphi_q &\varphi_d + \Phi & 0
    \end{bmatrix}\!,\, 
    R(x) \!=\! \begin{bmatrix}
        r & 0 & 0 \\
        0 & r & 0 \\
        0 & 0 & \beta
    \end{bmatrix}, \nonumber
\end{align}
with $x = \begin{bmatrix}
    \varphi_d &
    \varphi_q &
    p
\end{bmatrix}^\top$, $\varphi_d$ and $\varphi_q$ are magnetic fluxes while $p$ is the rotor momentum.
\begin{table}[h!]
\centering
\renewcommand{\arraystretch}{1.5}
\begin{tabular}{|c |c |c |c |}
\toprule
\hline
symbol& meaning & value & unit \\
\hline\hline
$J_m$ & inertia & 0.012 & $[{\rm kg \cdot m^2}]$ \\   $L$ & phase inductance & $3.8 \cdot 10^{-3}$ & $[{\rm H}]$ \\ 
$\beta$ & viscous friction coefficient & 0.0026 & $[{\rm Nms/rad}]$ \\ $r$ & phase resistance & 0.225 & $[\Omega]$ \\
$\Phi$ & constant rotor magnetic flux & 0.17 & $[{\rm Wb}]$ \\
\hline
\bottomrule
\end{tabular}
\caption{Permanent magnet synchronous motor system parameters}
\label{tab:PMSM}
\end{table}
The parameter values are shown in Table~\ref{tab:PMSM}.

\subsection{Wiener-Hammerstein system}
The Wiener-Hammerstein system is a widely studied block-oriented structure composed of two linear time-invariant~(LTI) blocks with a static nonlinearity in between. The interaction between the nonlinear element and the surrounding dynamic blocks makes system identification particularly challenging. The Wiener-Hammerstein system that we consider here is described in~\cite{SCHOUKENS2017446} and is affected primarily by process noise, which enters upstream of the static nonlinearity. Two additional but much less significant noise sources are also present in the input and output measurement channels. The input-output measurement data that is used for training is taken from~\cite{doi_wiener_hammer_data}. 

\subsection{Silverbox system}

The Silverbox system is an electronic implementation of the Duffing oscillator. It comprises a second-order LTI system with a cubic static nonlinearity embedded in a feedback loop~\cite{WigrenSchoukens13}. The input–output data considered here is obtained from~\cite{silverbox_data}.
\color{black}

\section{Choosing the number of breaks}
{We will discuss in this section how to choose the number of breaks $b$ when approximating broken NSLs. In general, there is no universally correct rule. A higher number of breaks is generally expected to decrease the interpolation error of the lower envelope $E_O$ that is expressed through the margin $M$ in~\eqref{eqn:margin}.
}

\begin{figure*}[htbp!] 
    \centering

    \begin{subfigure}[b]{0.48\textwidth}
        \centering
        \includegraphics[width=\linewidth]{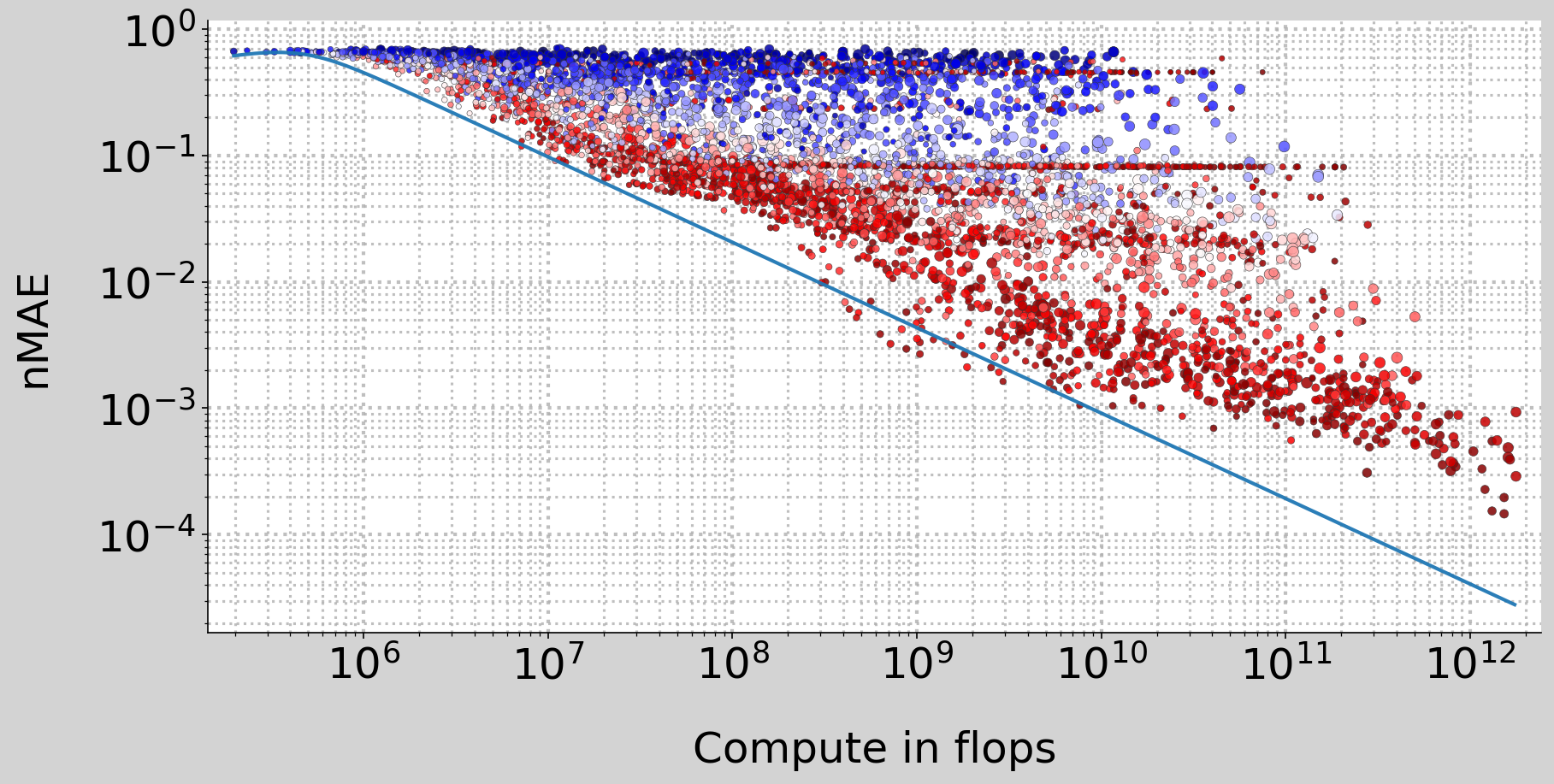}
        \caption{}
        \label{fig:ab-a}
    \end{subfigure}\hfill
    \begin{subfigure}[b]{0.48\textwidth}
        \centering
        \includegraphics[width=\linewidth]{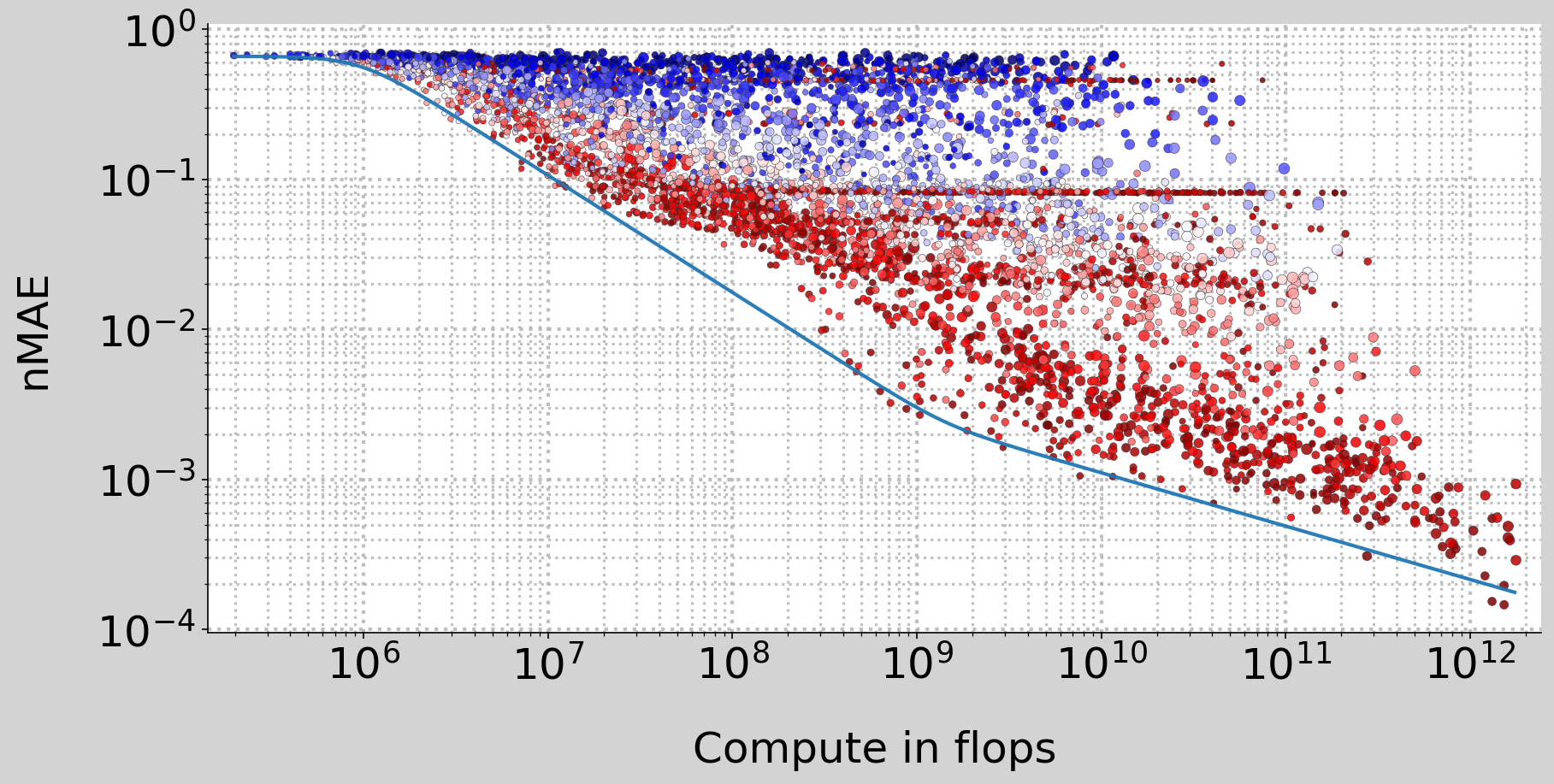}
        \caption{}
        \label{fig:ab-b}
    \end{subfigure}

    \vspace{0.5em}

    \begin{subfigure}[b]{0.48\textwidth}
        \centering
        \includegraphics[width=\linewidth]{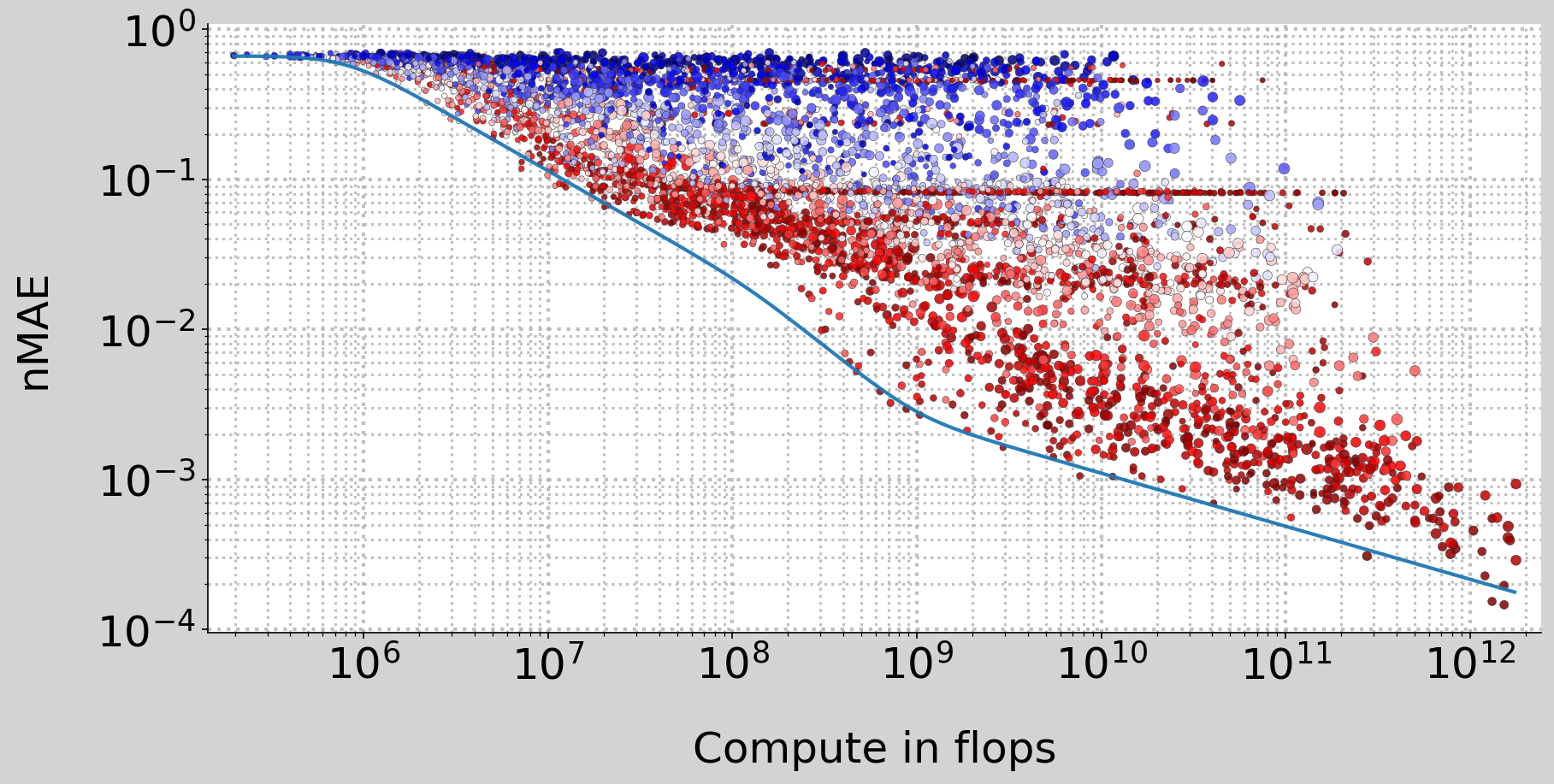}
        \caption{}
        \label{fig:ab-c}
    \end{subfigure}\hfill
    \begin{subfigure}[b]{0.48\textwidth}
        \centering
        \includegraphics[width=\linewidth]{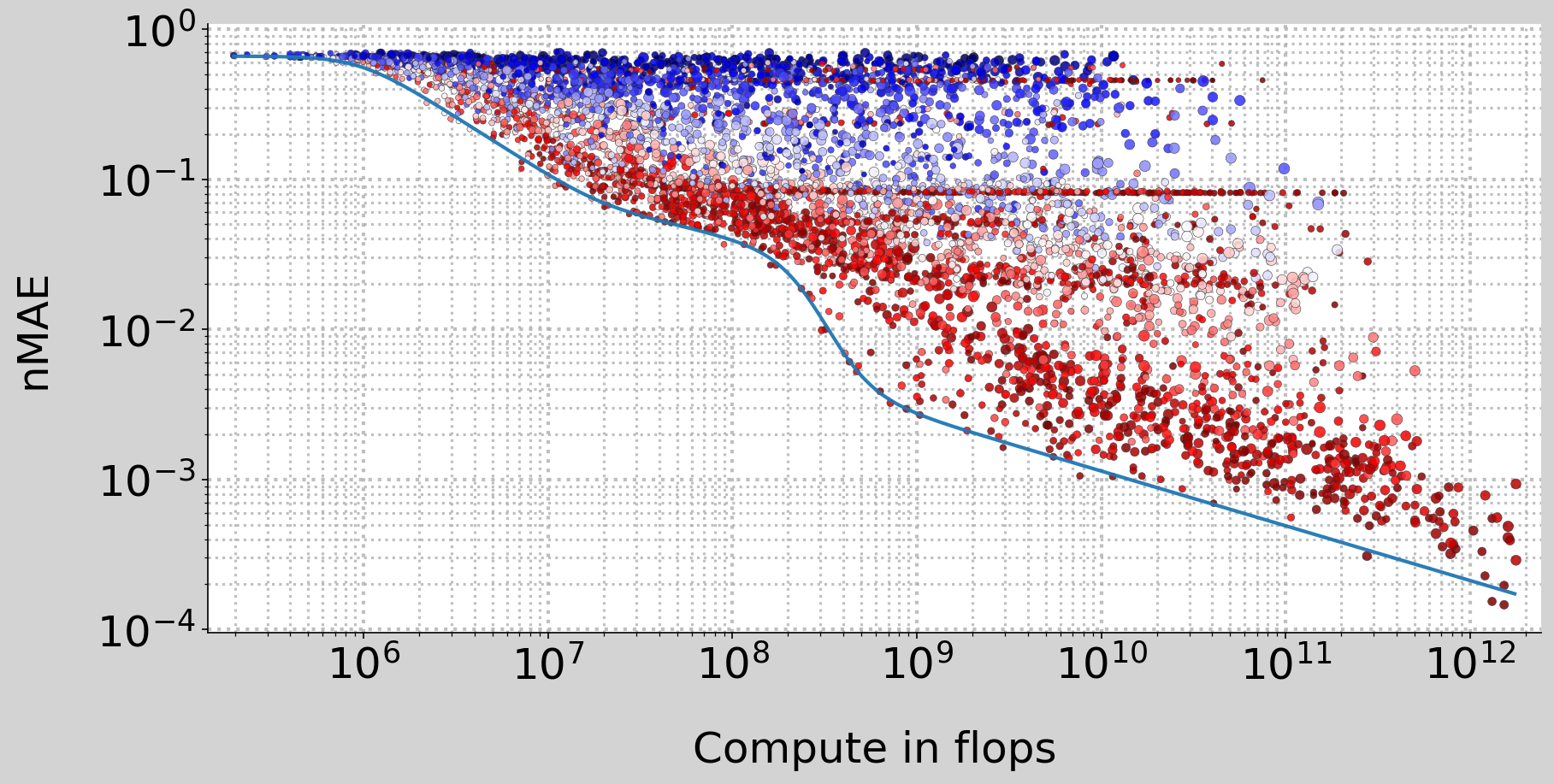}
        \caption{}
        \label{fig:ab-d}
    \end{subfigure}

    \caption{{Various compute-$\operatorname{nMAE}$ NSLs using input-state-output data for a varying number of breaks ranging from $b=1$ in (a) up to $b=4$ in (d).}}
    \label{fig:ab-margin}
\end{figure*}
{ Figure~\ref{fig:ab-margin} shows the compute-$\operatorname{nMAE}$ NSLs for the unstructured architecture on the spring system with different numbers of breaks.
Figure~\ref{fig:ab-a} shows the corresponding NSL with one break. We can clearly see that the number of breaks is not enough to accurately approximate the lower envelope $E_O$ of the set of outcomes $O$. In this case, the margin \eqref{eqn:margin} equals $M=0.68$.}
{Figure~\ref{fig:ab-b} shows the the corresponding NSL using two breaks.
In this case, we can see that the NSL fits the data more accurately which also corresponds to the lower margin $M=0.12$.
However, in the region of roughly $10^8\, \mathrm{flops}$ there is a gap between the point cloud and the blue curve representing the approximated NSL. This can be interpreted as acceptable if one assumes that this region can be considered an outlier-region.
}
{Figure~\ref{fig:ab-c} shows the the corresponding NSL using three breaks. In this case, we can see that adding one more break does not yield a significant improvement and only lowers the margin $M$ from $M=0.12$ to $M=0.08$.}
{Figure~\ref{fig:ab-d} shows the corresponding NSL using four breaks. In this case, we can see that the gap observed earlier does not exist anymore and that the blue curve fits the lower envelope $E_O$ almost perfectly. This can also be verified by considering the significantly lower margin $M=0.03$.}

{Summarizing this discussion, we have seen that there is a trade-off between using more breaks (more accurate approximation of the lower envelope and lower margin) and the simplicity of an NSL. Moreover, the right choice of the number of breaks can become ambiguous depending on which regions are considered to be 'outlier-regions'.}

\section{Overview of neural scaling laws using input-state-output data}\label{app:ov}

\begin{figure*}[t] 
    \centering

    \begin{minipage}[b]{0.32\textwidth}
        \centering
        \includegraphics[width=\linewidth]{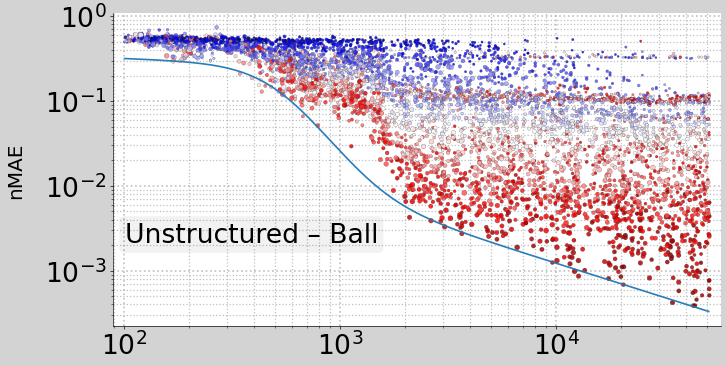}
    \end{minipage}\hfill
    \begin{minipage}[b]{0.32\textwidth}
        \centering
        \includegraphics[width=\linewidth]{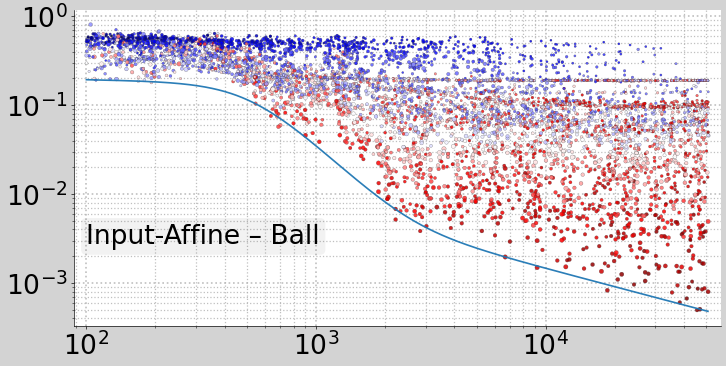}
    \end{minipage}\hfill
    \begin{minipage}[b]{0.32\textwidth}
        \centering
        \includegraphics[width=\linewidth]{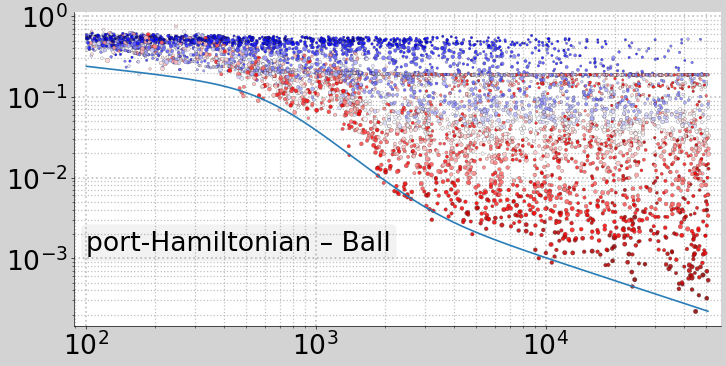}
    \end{minipage}

    \begin{minipage}[b]{0.32\textwidth}
        \centering
        \includegraphics[width=\linewidth]{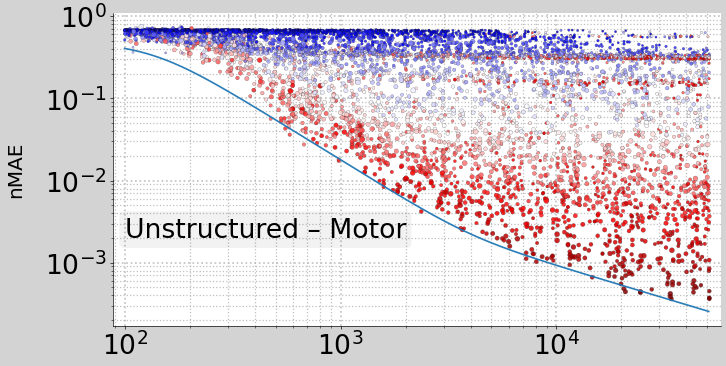}
    \end{minipage}\hfill
    \begin{minipage}[b]{0.32\textwidth}
        \centering
        \includegraphics[width=\linewidth]{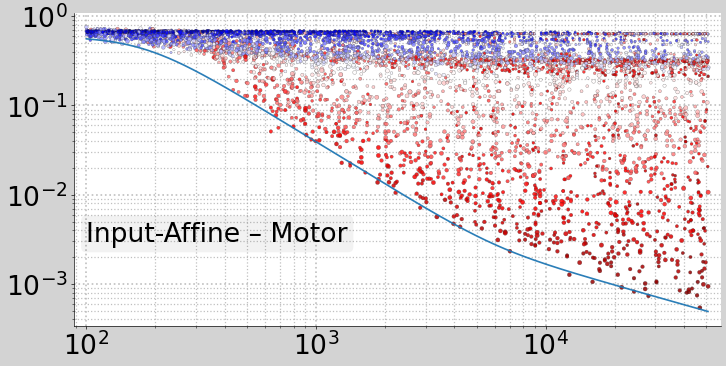}
    \end{minipage}\hfill
    \begin{minipage}[b]{0.32\textwidth}
        \centering
        \includegraphics[width=\linewidth]{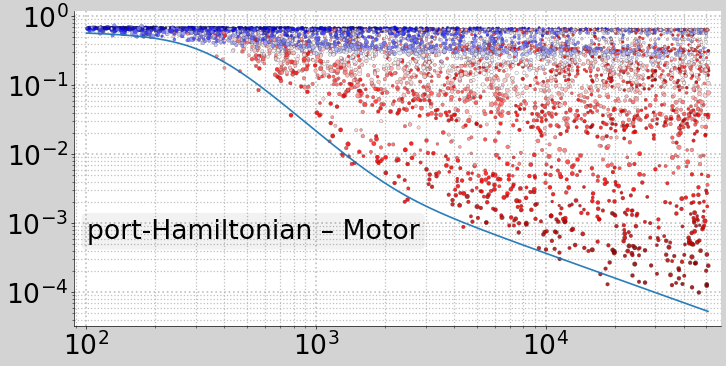}
    \end{minipage}

    \begin{minipage}[b]{0.32\textwidth}
        \centering
        \includegraphics[width=\linewidth]{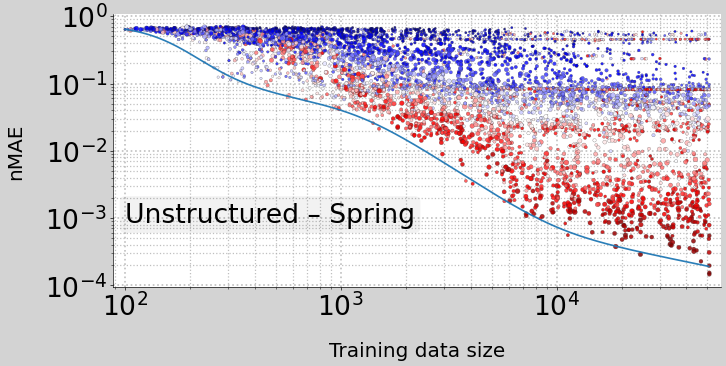}
    \end{minipage}\hfill
    \begin{minipage}[b]{0.32\textwidth}
        \centering
        \includegraphics[width=\linewidth]{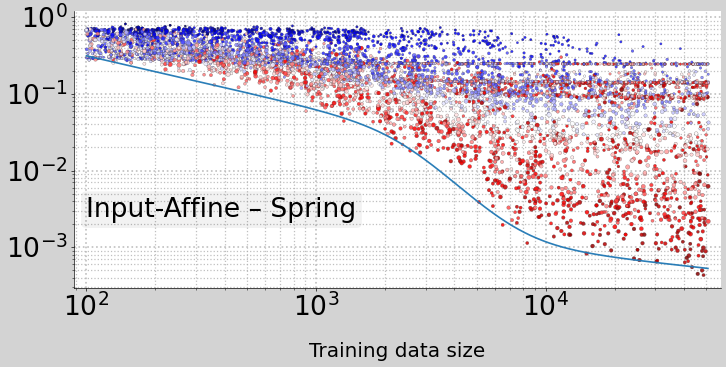}
    \end{minipage}\hfill
    \begin{minipage}[b]{0.32\textwidth}
        \centering
        \includegraphics[width=\linewidth]{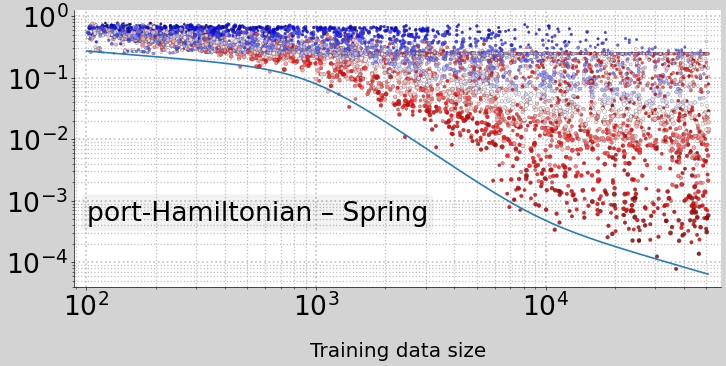}
    \end{minipage}

    \caption{All data-$\operatorname{nMAE}$ NSLs {using input-state-output data}}
    \label{fig:data-nmae}
\end{figure*}

\begin{figure*}[t] 
    \centering

    \begin{minipage}[b]{0.32\textwidth}
        \centering
        \includegraphics[width=\linewidth]{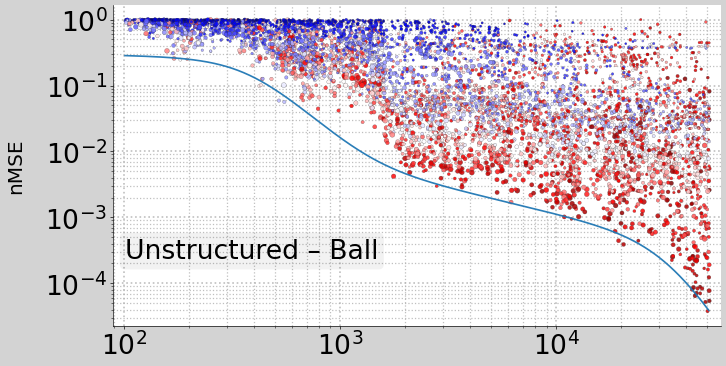}
    \end{minipage}\hfill
    \begin{minipage}[b]{0.32\textwidth}
        \centering
        \includegraphics[width=\linewidth]{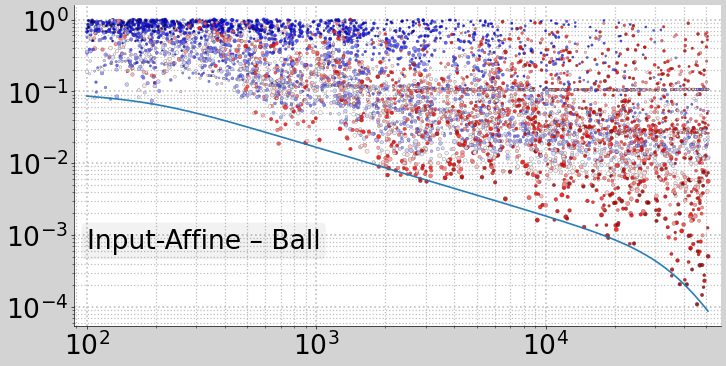}
    \end{minipage}\hfill
    \begin{minipage}[b]{0.32\textwidth}
        \centering
        \includegraphics[width=\linewidth]{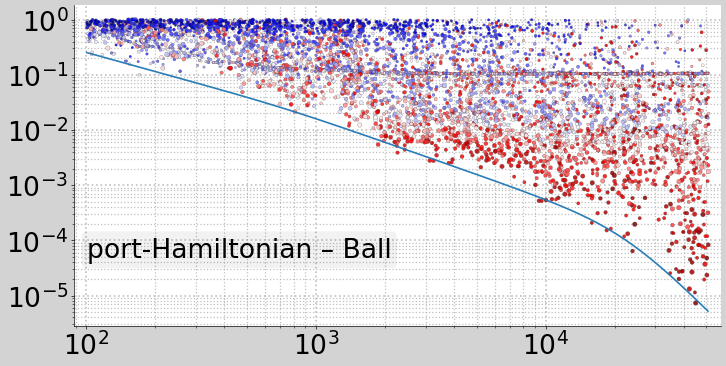}
    \end{minipage}

    \begin{minipage}[b]{0.32\textwidth}
        \centering
        \includegraphics[width=\linewidth]{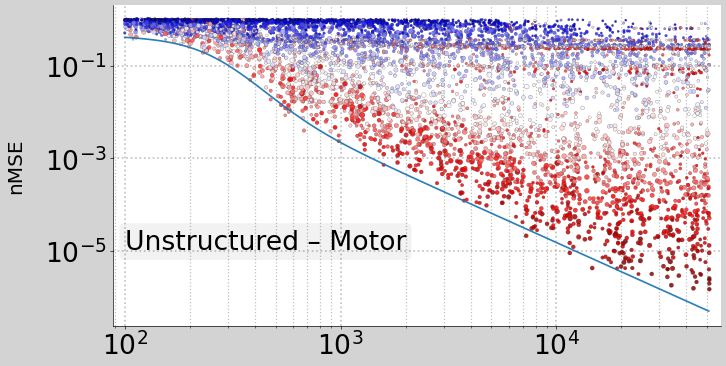}
    \end{minipage}\hfill
    \begin{minipage}[b]{0.32\textwidth}
        \centering
        \includegraphics[width=\linewidth]{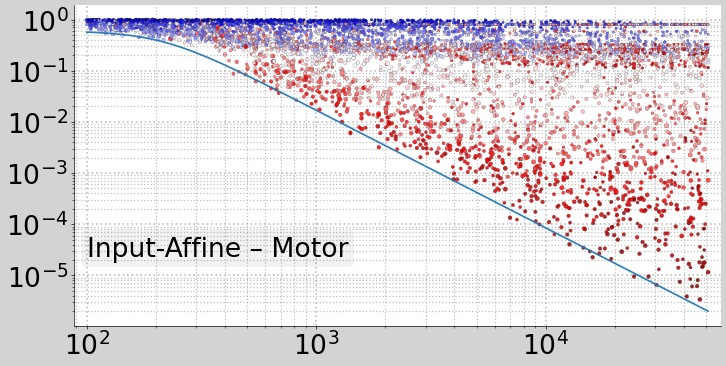}
    \end{minipage}\hfill
    \begin{minipage}[b]{0.32\textwidth}
        \centering
        \includegraphics[width=\linewidth]{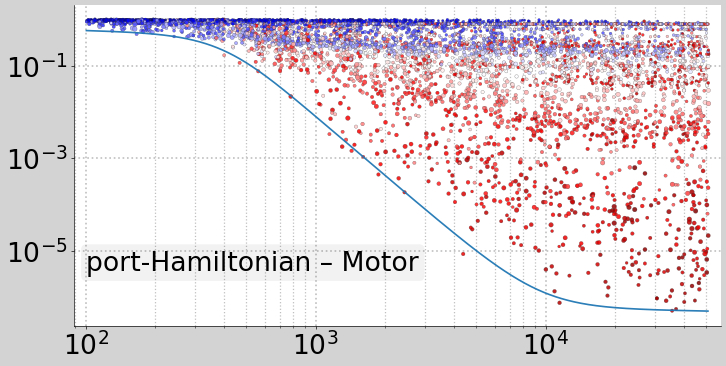}
    \end{minipage}

    \begin{minipage}[b]{0.32\textwidth}
        \centering
        \includegraphics[width=\linewidth]{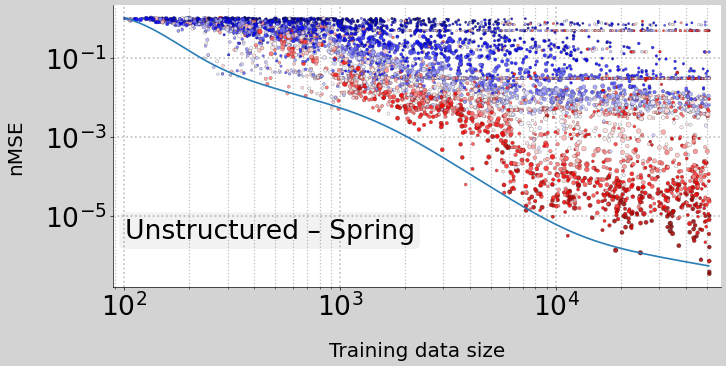}
    \end{minipage}\hfill
    \begin{minipage}[b]{0.32\textwidth}
        \centering
        \includegraphics[width=\linewidth]{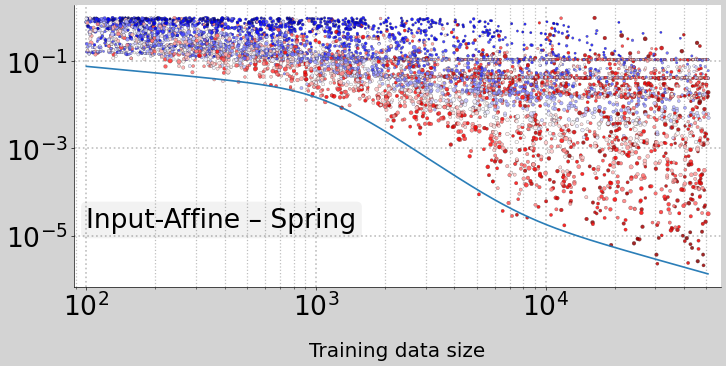}
    \end{minipage}\hfill
    \begin{minipage}[b]{0.32\textwidth}
        \centering
        \includegraphics[width=\linewidth]{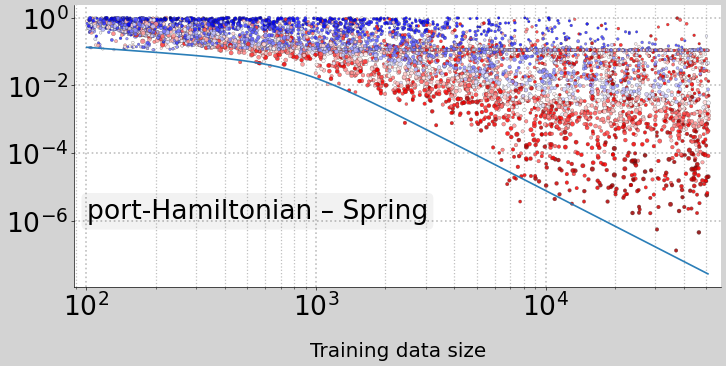}
    \end{minipage}

    \caption{All data-$\operatorname{nMSE}$ NSLs using {using input-state-output data}}
    \label{fig:data-nmse}
\end{figure*}

\begin{figure*}[t] 
    \centering

    \begin{minipage}[b]{0.32\textwidth}
        \centering
        \includegraphics[width=\linewidth]{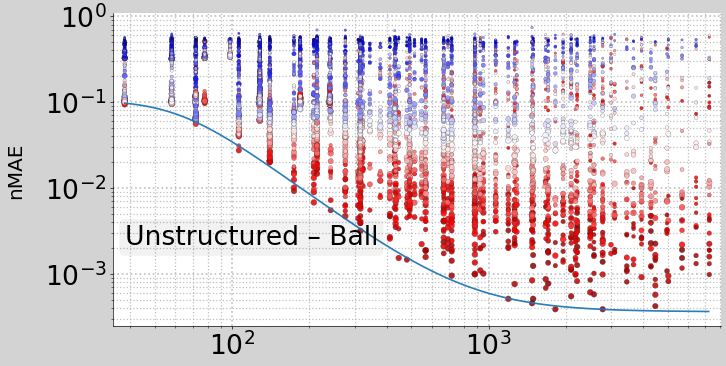}
    \end{minipage}\hfill
    \begin{minipage}[b]{0.32\textwidth}
        \centering
        \includegraphics[width=\linewidth]{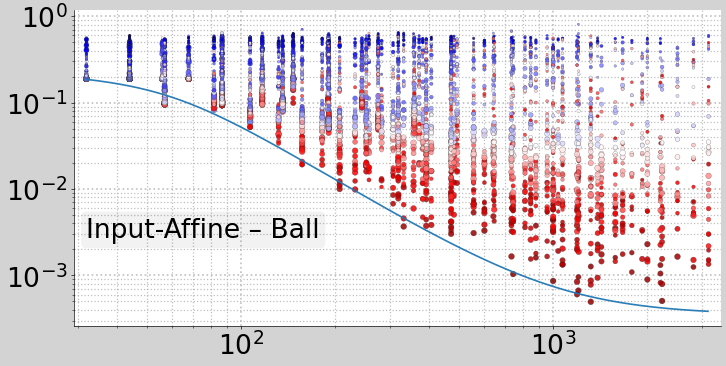}
    \end{minipage}\hfill
    \begin{minipage}[b]{0.32\textwidth}
        \centering
        \includegraphics[width=\linewidth]{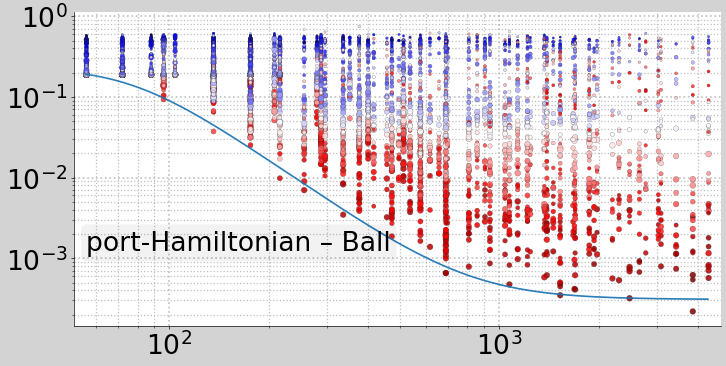}
    \end{minipage}

    \begin{minipage}[b]{0.32\textwidth}
        \centering
        \includegraphics[width=\linewidth]{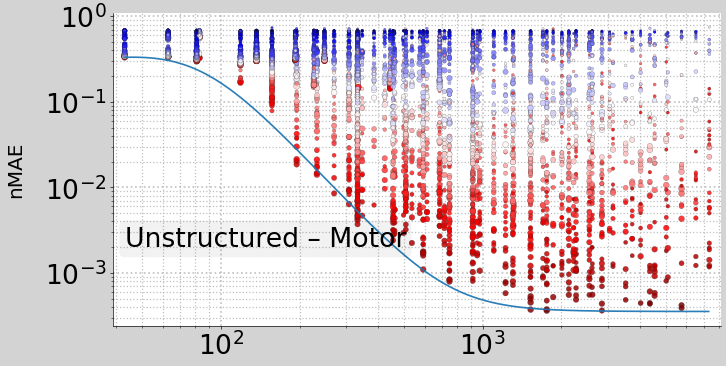}
    \end{minipage}\hfill
    \begin{minipage}[b]{0.32\textwidth}
        \centering
        \includegraphics[width=\linewidth]{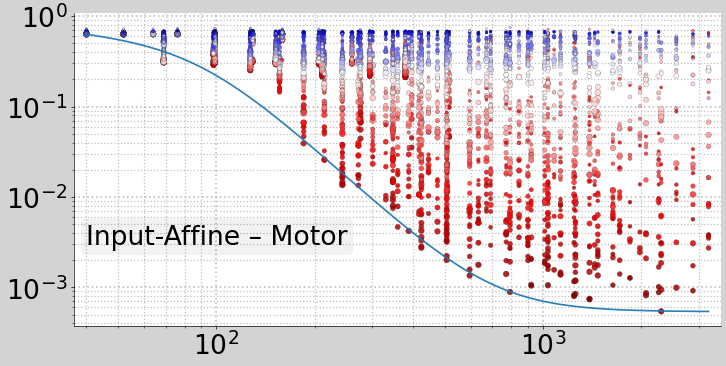}
    \end{minipage}\hfill
    \begin{minipage}[b]{0.32\textwidth}
        \centering
        \includegraphics[width=\linewidth]{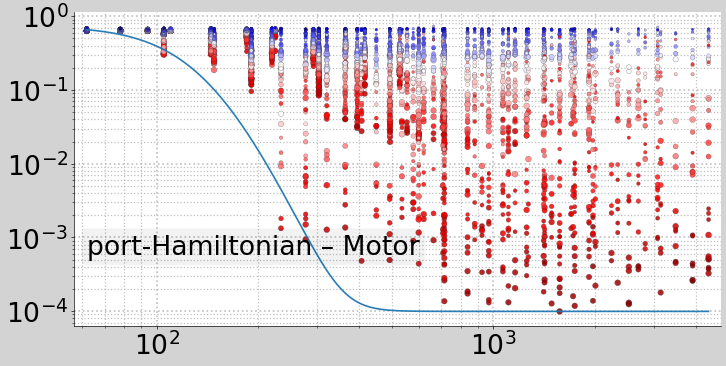}
    \end{minipage}

    \begin{minipage}[b]{0.32\textwidth}
        \centering
        \includegraphics[width=\linewidth]{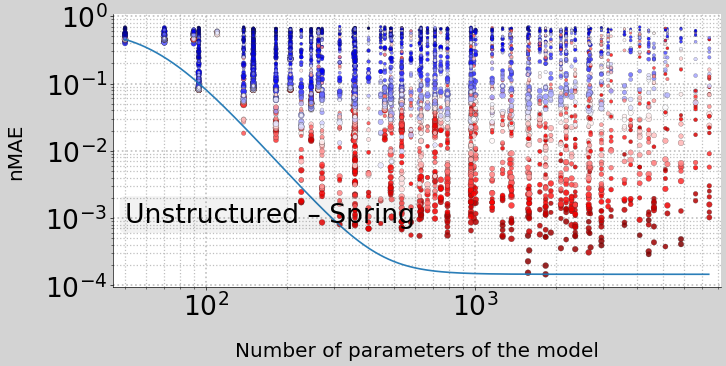}
    \end{minipage}\hfill
    \begin{minipage}[b]{0.32\textwidth}
        \centering
        \includegraphics[width=\linewidth]{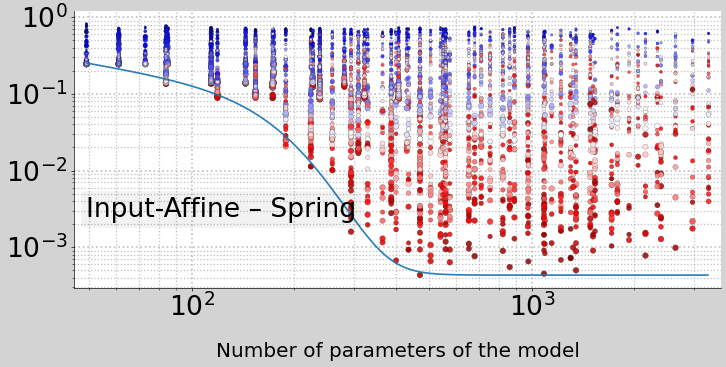}
    \end{minipage}\hfill
    \begin{minipage}[b]{0.32\textwidth}
        \centering
        \includegraphics[width=\linewidth]{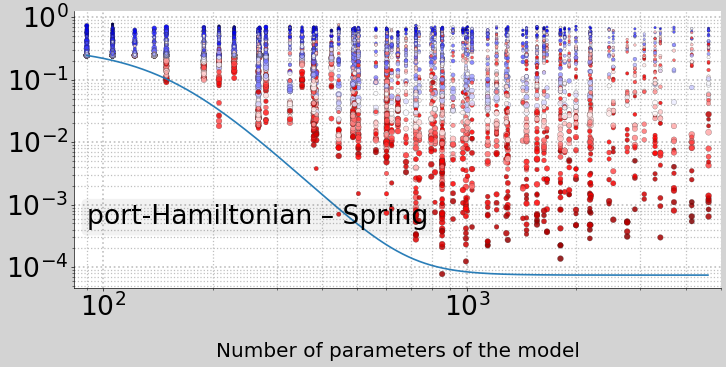}
    \end{minipage}

    \caption{All model-$\operatorname{nMAE}$ NSLs {using input-state-output data}}
    \label{fig:model-nmae}
\end{figure*}

\begin{figure*}[t] 
    \centering

    \begin{minipage}[b]{0.32\textwidth}
        \centering
        \includegraphics[width=\linewidth]{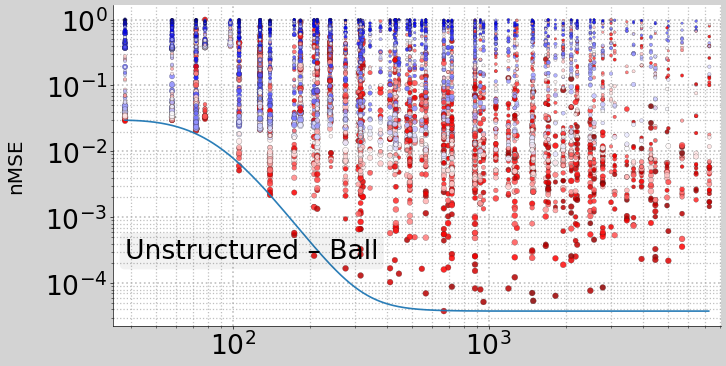}
    \end{minipage}\hfill
    \begin{minipage}[b]{0.32\textwidth}
        \centering
        \includegraphics[width=\linewidth]{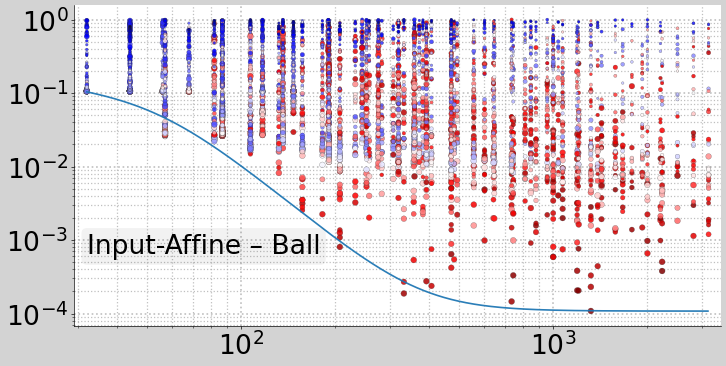}
    \end{minipage}\hfill
    \begin{minipage}[b]{0.32\textwidth}
        \centering
        \includegraphics[width=\linewidth]{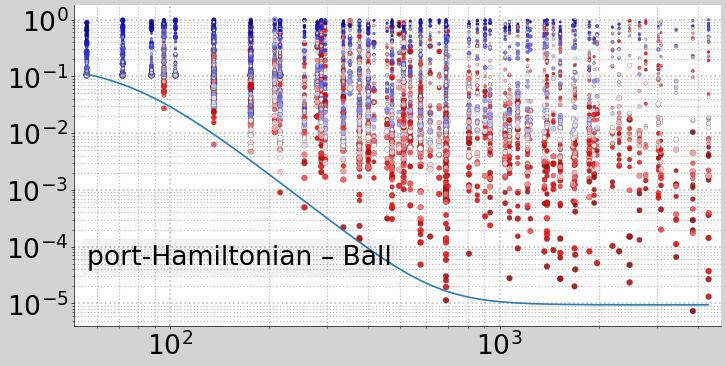}
    \end{minipage}

    \begin{minipage}[b]{0.32\textwidth}
        \centering
        \includegraphics[width=\linewidth]{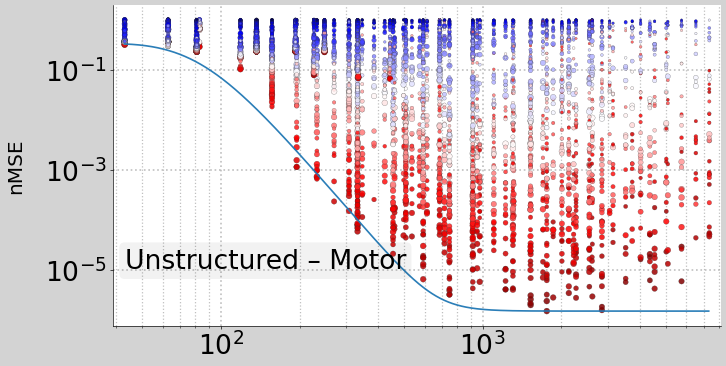}
    \end{minipage}\hfill
    \begin{minipage}[b]{0.32\textwidth}
        \centering
        \includegraphics[width=\linewidth]{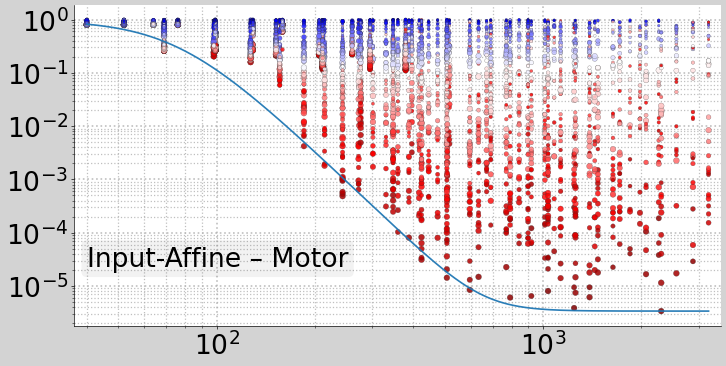}
    \end{minipage}\hfill
    \begin{minipage}[b]{0.32\textwidth}
        \centering
        \includegraphics[width=\linewidth]{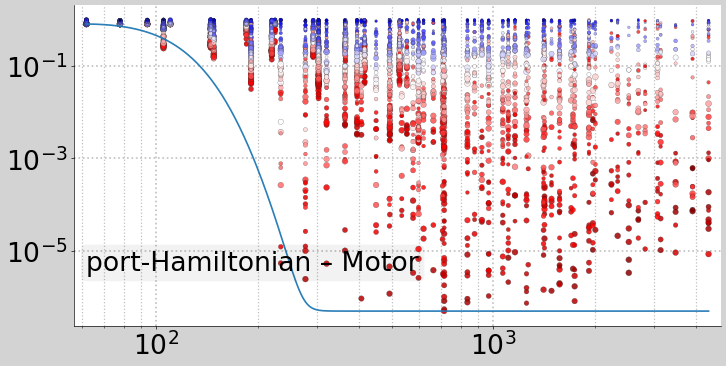}
    \end{minipage}

    \begin{minipage}[b]{0.32\textwidth}
        \centering
        \includegraphics[width=\linewidth]{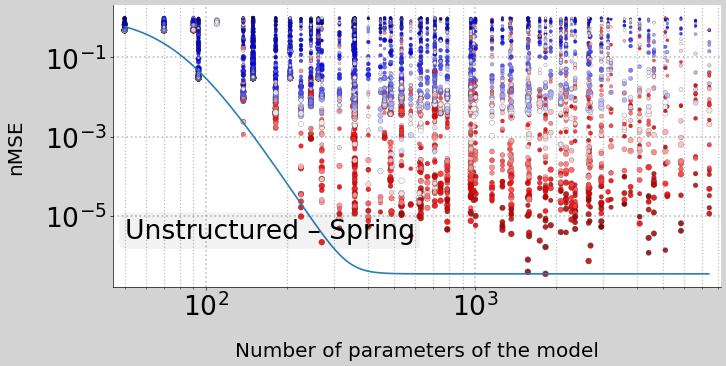}
    \end{minipage}\hfill
    \begin{minipage}[b]{0.32\textwidth}
        \centering
        \includegraphics[width=\linewidth]{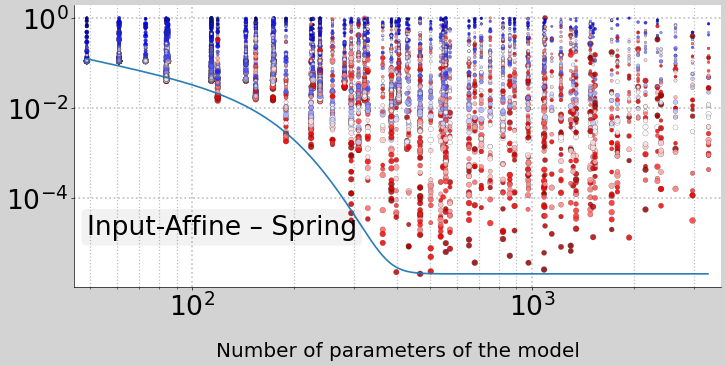}
    \end{minipage}\hfill
    \begin{minipage}[b]{0.32\textwidth}
        \centering
        \includegraphics[width=\linewidth]{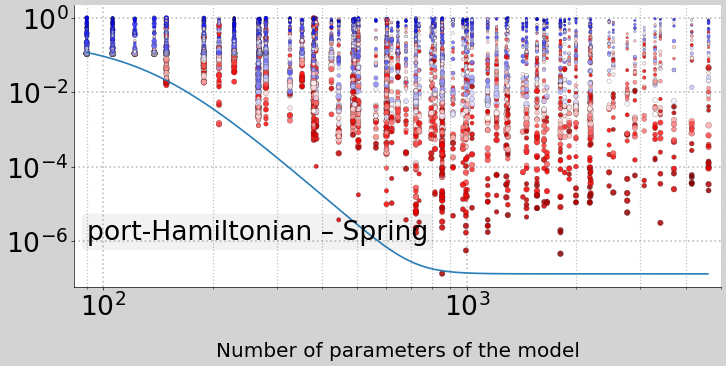}
    \end{minipage}

    \caption{All model-$\operatorname{nMSE}$ NSLs {using input-state-output data}}
    \label{fig:model-nmse}
\end{figure*}

In this section, we collect the figures of all considered NSLs {using input-state-output data}.
The compute-$\operatorname{nMAE}$ NSLs are shown in 
Figure~\ref{fig:comp-nmae}. All compute-$\operatorname{nMSE}$ NSLs are shown in Figure~\ref{fig:comp-nmse}. 
{We can see that these NSLs share a high amount of structural similarity. They all exhibit a large regimes of smooth scaling with more compute. Moreover, all of these NSLs need a minimum amount of compute to see any improvement beyond the initial plateau. In some cases, however, they differ in the number of breaks one can see. The unstructured architecture on the spring system for example has a relatively high number of $4$ breaks.
Figure~\ref{fig:comp-nmse} shows the
corresponding NSLs for the nMSE metric. In this case we can see that the observation made in Section~\ref{sec:compute} that using the nMSE instead of nMAE results in the high compute region having less point density is true.
}
The data-$\operatorname{nMAE}$ NSLs are shown in Figure~\ref{fig:data-nmae}. 
{Compared to the compute-nMAE NSLs we see more variation in this case. It should be noted, however that these NSLs often look similar for all the different architectures on a given dataset. On the ball system for example all three NSLs first show a phase without scaling, then a short window where fast progression appears and then enter a final regime of smooth scaling. This is in contrast to Figure~\ref{fig:data-nmse} where we see more variation in the data-nMAE NSLs compared to the corresponding the NSLs using the nMAE metric. Thus, the nMAE seems to be a more natural metric for NSLs than the nMSE.}
The model-$\operatorname{nMAE}$ NSLs are shown in Figure~\ref{fig:model-nmae}. All model-$\operatorname{nMSE}$ NSLs are shown in Figure~\ref{fig:model-nmse}.
{We can see in both of these figures that a certain minimum of parameters is necessary to obtain the optimal performance, but all of the NSLs saturate at relatively low parameter counts.}
{Note that we were able to find a neural scaling law for all model structures and all examples, as shown in the figures.} 
{
While a comparison between different architectures is not the main focus of this work, we note that NSLs enable us to get insights into the scaling behavior for different architectures on a given task.
For example if one considers the NSLs in Figure~\ref{fig:data-nmae}, the order of magnitude in accuracy for a given example between the unstructured, input-affine and port-Hamiltonian are relatively similar at the same data size. In some cases we get one order of magnitude better accuracy but this is not consistent for all examples.
}

\section{Overview of neural scaling laws using input-output data}\label{app:ov-io}
\begin{figure*}[htbp!] 
    \centering

    \begin{minipage}[b]{0.32\textwidth}
        \centering        \includegraphics[width=\linewidth]{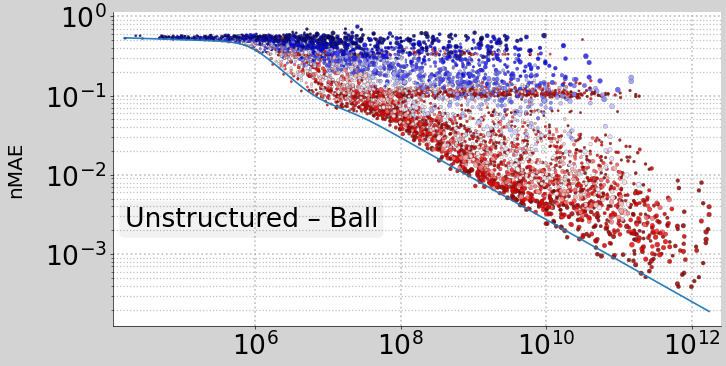}
    \end{minipage}\hfill
    \begin{minipage}[b]{0.32\textwidth}
        \centering
        \includegraphics[width=\linewidth]{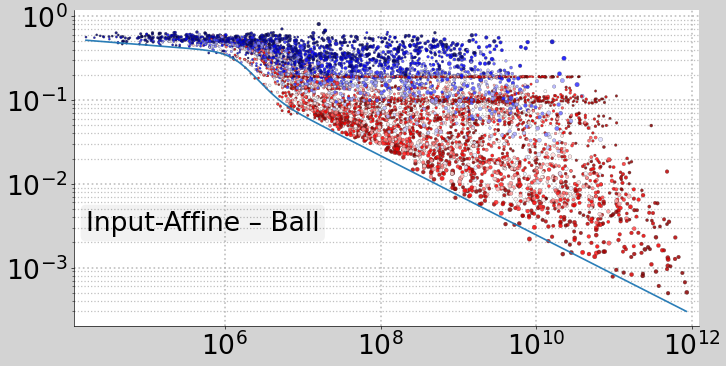}
    \end{minipage}\hfill
    \begin{minipage}[b]{0.32\textwidth}
        \centering
        \includegraphics[width=\linewidth]{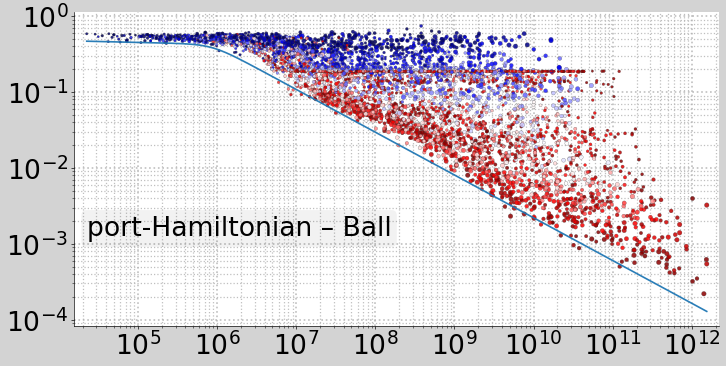}
    \end{minipage}

    \begin{minipage}[b]{0.32\textwidth}
        \centering
        \includegraphics[width=\linewidth]{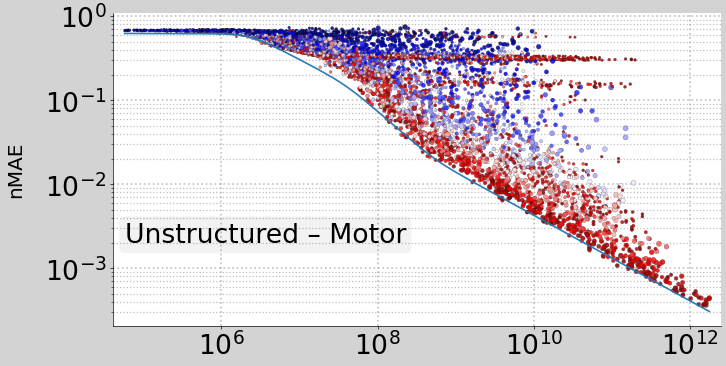}
    \end{minipage}\hfill
    \begin{minipage}[b]{0.32\textwidth}
        \centering
        \includegraphics[width=\linewidth]{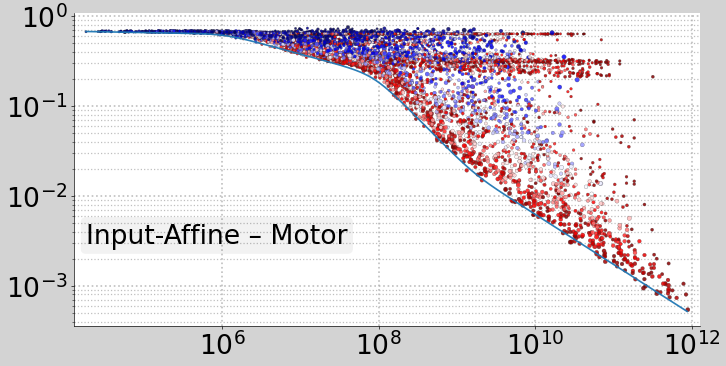}
    \end{minipage}\hfill
    \begin{minipage}[b]{0.32\textwidth}
        \centering
        \includegraphics[width=\linewidth]{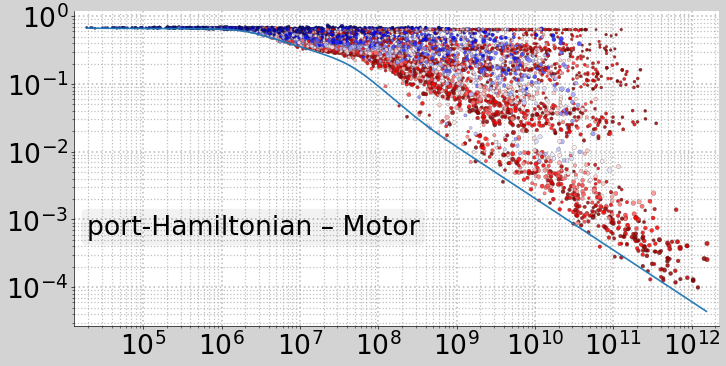}
    \end{minipage}

    \begin{minipage}[b]{0.32\textwidth}
        \centering
        \includegraphics[width=\linewidth]{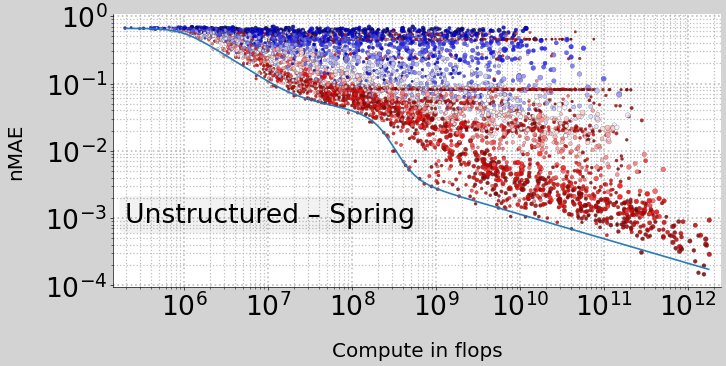}
    \end{minipage}\hfill
    \begin{minipage}[b]{0.32\textwidth}
        \centering
        \includegraphics[width=\linewidth]{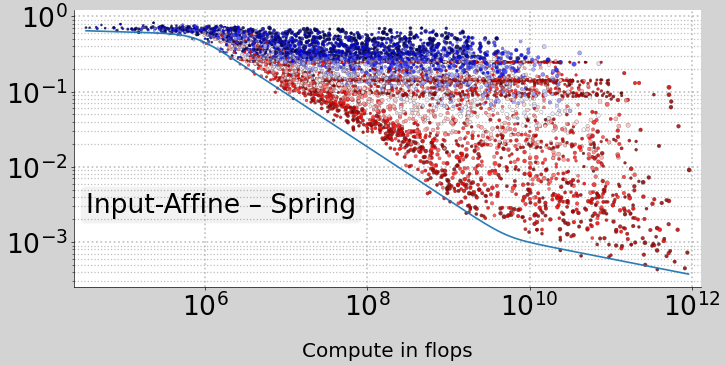}
    \end{minipage}\hfill
    \begin{minipage}[b]{0.32\textwidth}
        \centering
        \includegraphics[width=\linewidth]{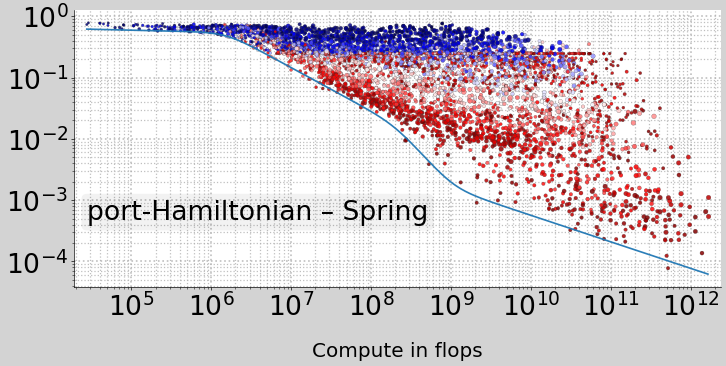}
    \end{minipage}
    
    \caption{All compute-$\operatorname{nMAE}$ NSLs {using input-state-output data}}
    \label{fig:comp-nmae}
\end{figure*}

\begin{figure*}[htbp!] 
    \centering

    \begin{minipage}[b]{0.32\textwidth}
        \centering
        \includegraphics[width=\linewidth]{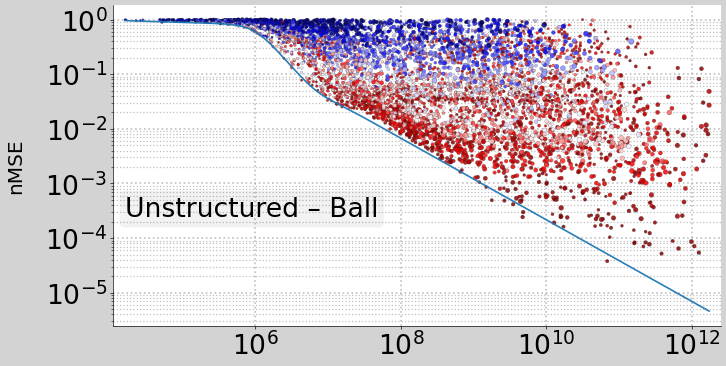}
    \end{minipage}\hfill
    \begin{minipage}[b]{0.32\textwidth}
        \centering
        \includegraphics[width=\linewidth]{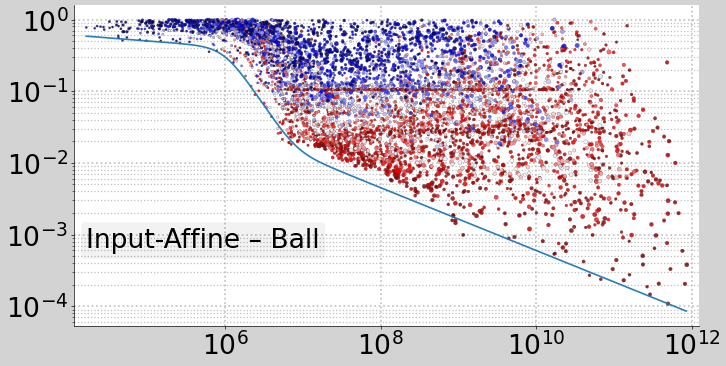}
    \end{minipage}\hfill
    \begin{minipage}[b]{0.32\textwidth}
        \centering
        \includegraphics[width=\linewidth]{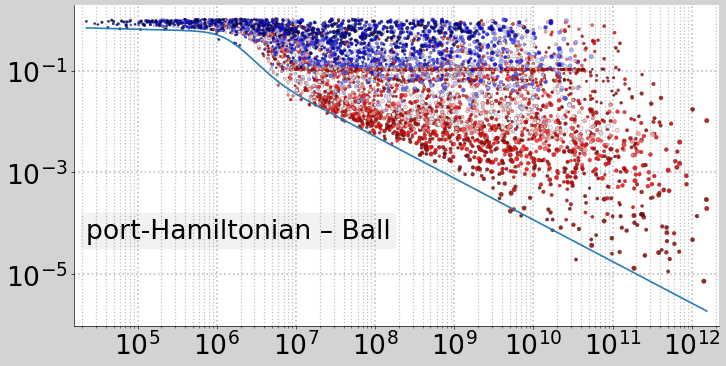}
    \end{minipage}

    \begin{minipage}[b]{0.32\textwidth}
        \centering
        \includegraphics[width=\linewidth]{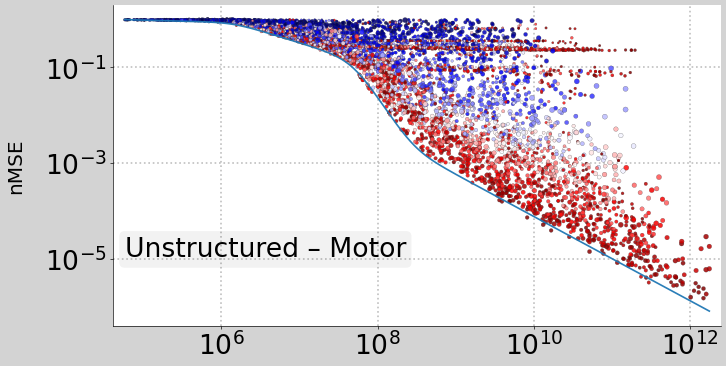}
    \end{minipage}\hfill
    \begin{minipage}[b]{0.32\textwidth}
        \centering
        \includegraphics[width=\linewidth]{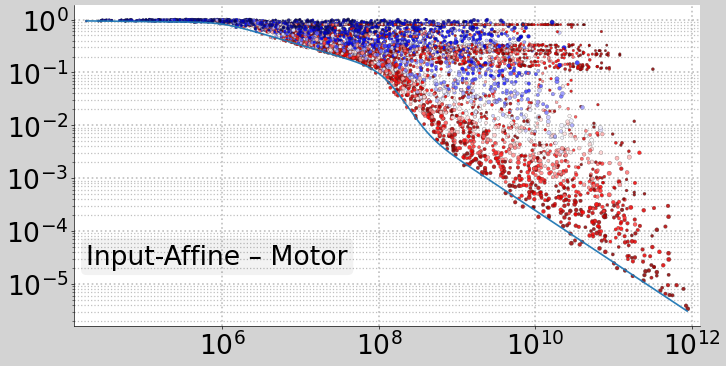}
    \end{minipage}\hfill
    \begin{minipage}[b]{0.32\textwidth}
        \centering
        \includegraphics[width=\linewidth]{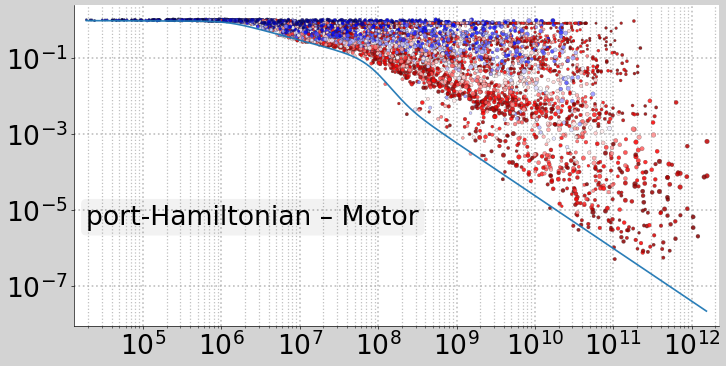}
    \end{minipage}

    \begin{minipage}[b]{0.32\textwidth}
        \centering
        \includegraphics[width=\linewidth]{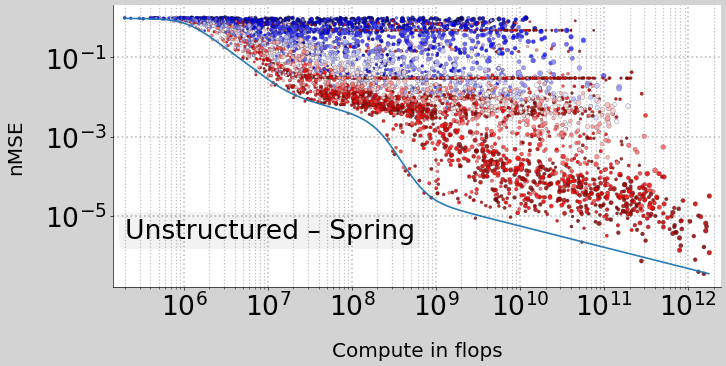}
    \end{minipage}\hfill
    \begin{minipage}[b]{0.32\textwidth}
        \centering
        \includegraphics[width=\linewidth]{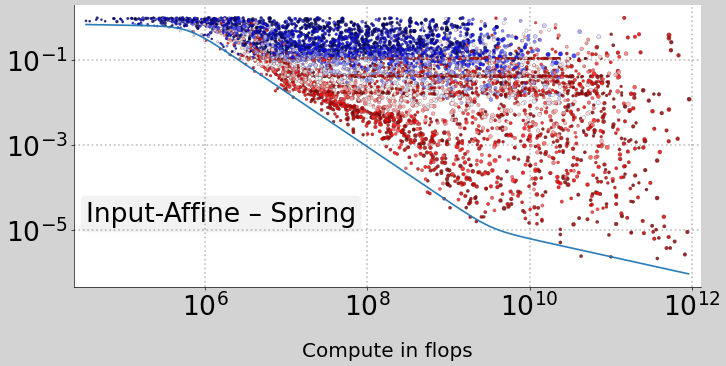}
    \end{minipage}\hfill
    \begin{minipage}[b]{0.32\textwidth}
        \centering
        \includegraphics[width=\linewidth]{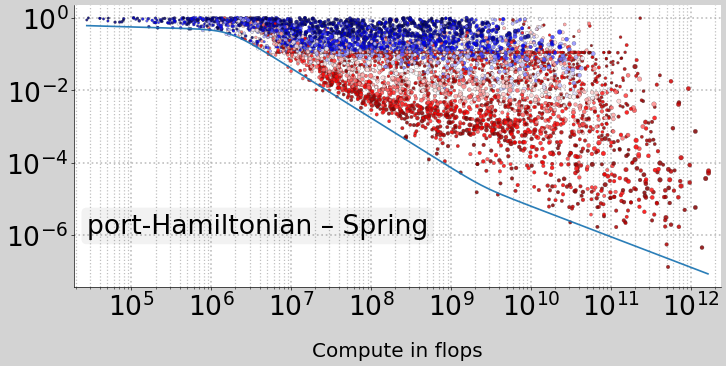}
    \end{minipage}

    \caption{All compute-    $\operatorname{nMSE}$ NSLs {using input-state-output data}}
    \label{fig:comp-nmse}
\end{figure*}

\begin{figure*}[htbp!] 
    \centering

    \begin{minipage}[b]{0.48\textwidth}
        \centering        \includegraphics[width=\linewidth]{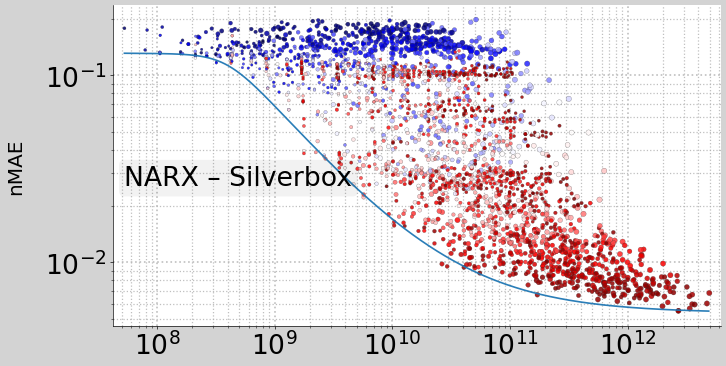}
    \end{minipage}\hfill
    \begin{minipage}[b]{0.48\textwidth}
        \centering
        \includegraphics[width=\linewidth]{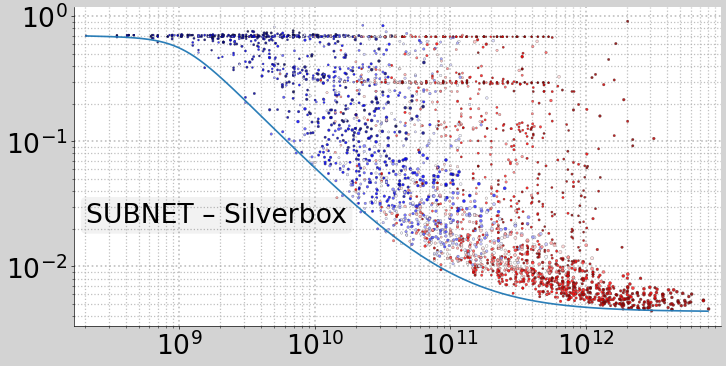}
    \end{minipage}

    \begin{minipage}[b]{0.48\textwidth}
        \centering
        \includegraphics[width=\linewidth]{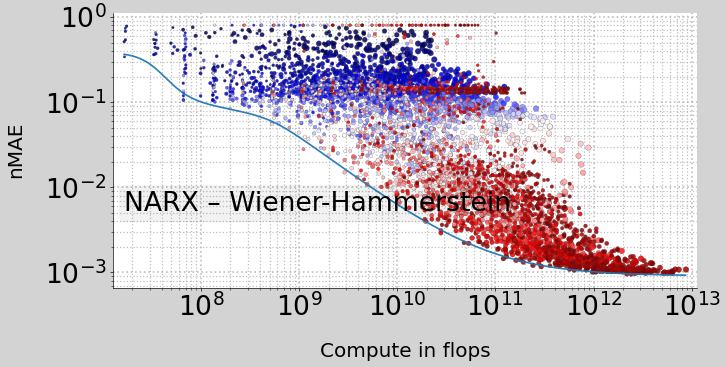}
    \end{minipage}\hfill
    \begin{minipage}[b]{0.48\textwidth}
        \centering
        \includegraphics[width=\linewidth]{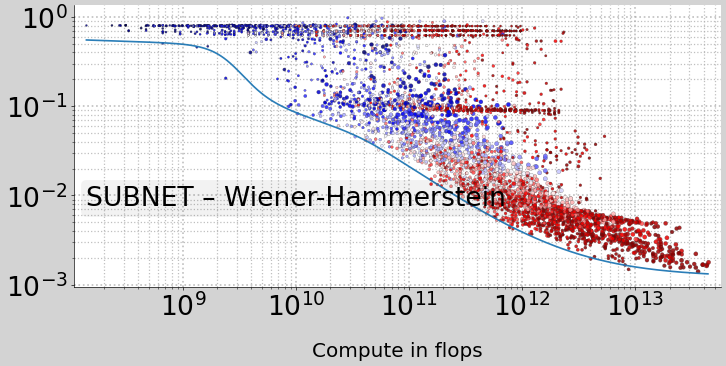}
    \end{minipage}

    \caption{{All compute-$\operatorname{nMAE}$ NSLs using input-output data.}}
    \label{fig:comp-io-nmae}
\end{figure*}

\begin{figure*}[htbp!] 
    \centering

    \begin{minipage}[b]{0.48\textwidth}
        \centering        \includegraphics[width=\linewidth]{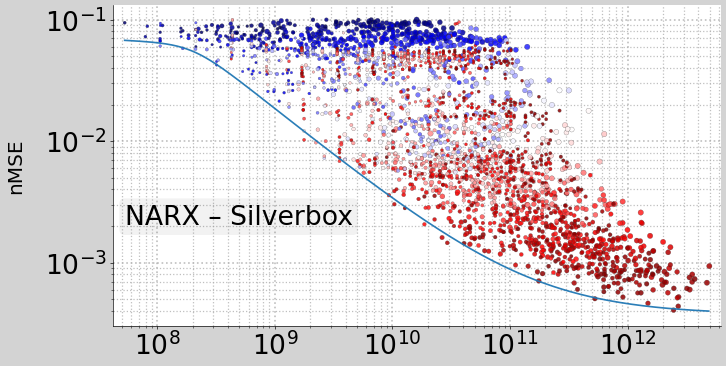}
    \end{minipage}\hfill
    \begin{minipage}[b]{0.48\textwidth}
        \centering
        \includegraphics[width=\linewidth]{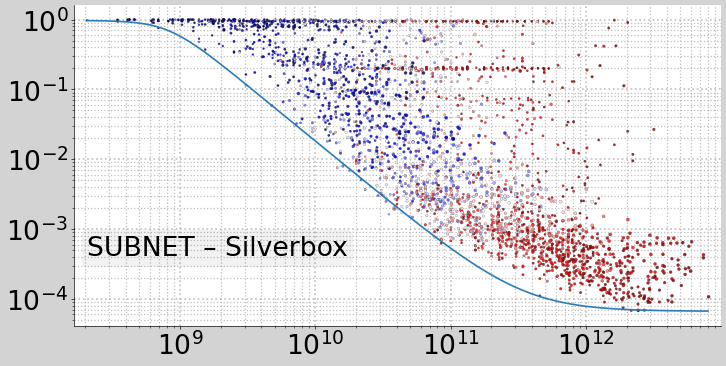}
    \end{minipage}

    \begin{minipage}[b]{0.48\textwidth}
        \centering
        \includegraphics[width=\linewidth]{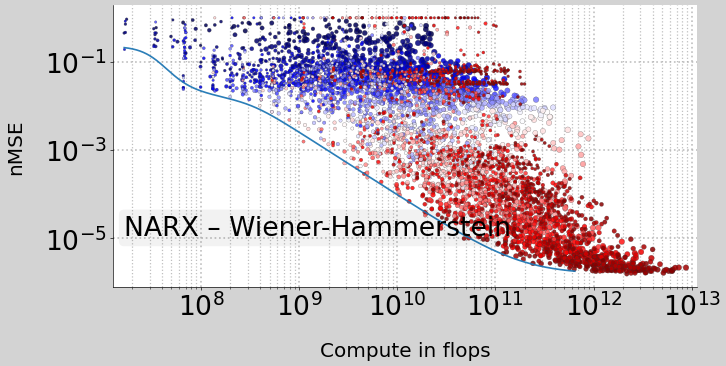}
    \end{minipage}\hfill
    \begin{minipage}[b]{0.48\textwidth}
        \centering
        \includegraphics[width=\linewidth]{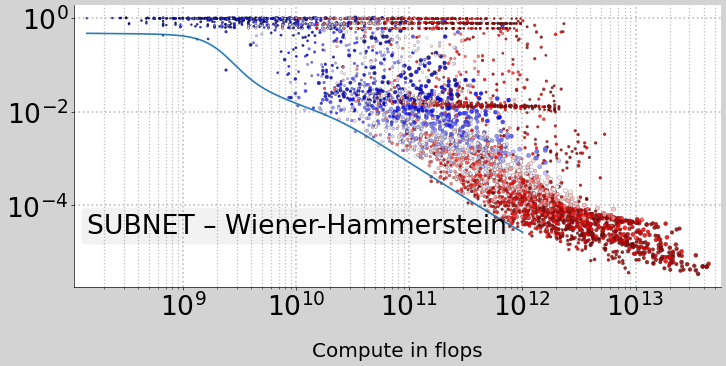}
    \end{minipage}

    \caption{{All compute-$\operatorname{nMSE}$ NSLs using input-output data.}}
    \label{fig:comp-io-nmse}
\end{figure*}

\begin{figure*}[htbp!] 
    \centering

    \begin{minipage}[b]{0.48\textwidth}
        \centering        \includegraphics[width=\linewidth]{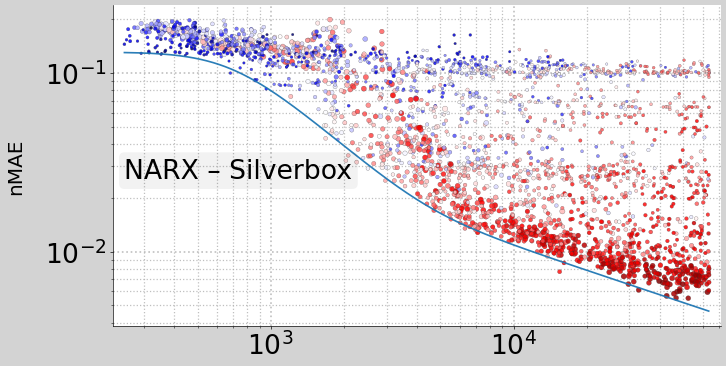}
    \end{minipage}\hfill
    \begin{minipage}[b]{0.48\textwidth}
        \centering
        \includegraphics[width=\linewidth]{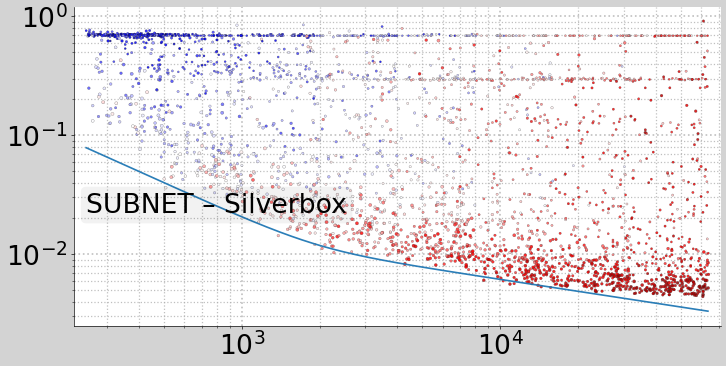}
    \end{minipage}

    \begin{minipage}[b]{0.48\textwidth}
        \centering
        \includegraphics[width=\linewidth]{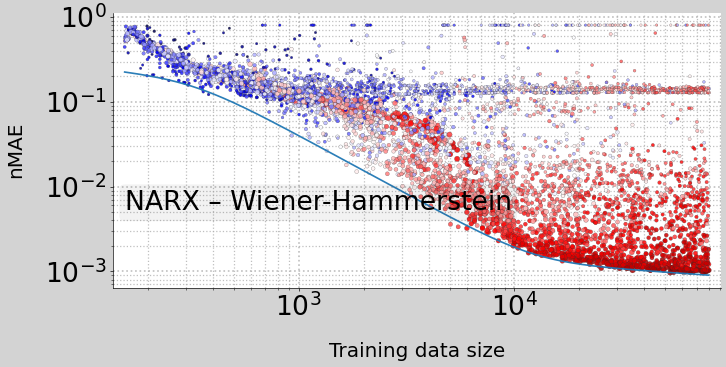}
    \end{minipage}\hfill
    \begin{minipage}[b]{0.48\textwidth}
        \centering
        \includegraphics[width=\linewidth]{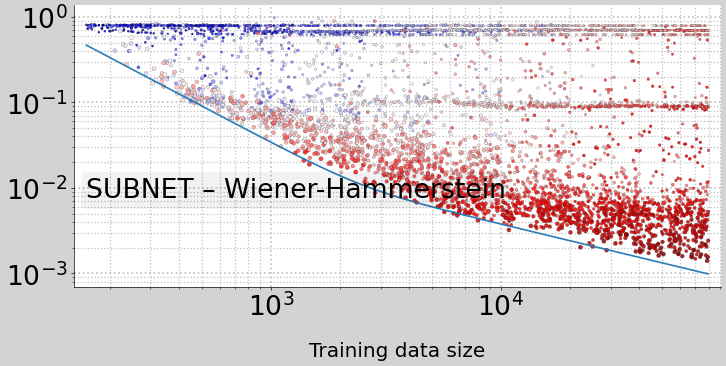}
    \end{minipage}

    \caption{{All data-$\operatorname{nMAE}$ NSLs using input-output data.}}
    \label{fig:data-io-nmae}
\end{figure*}

\begin{figure*}[htbp!] 
    \centering

    \begin{minipage}[b]{0.48\textwidth}
        \centering        \includegraphics[width=\linewidth]{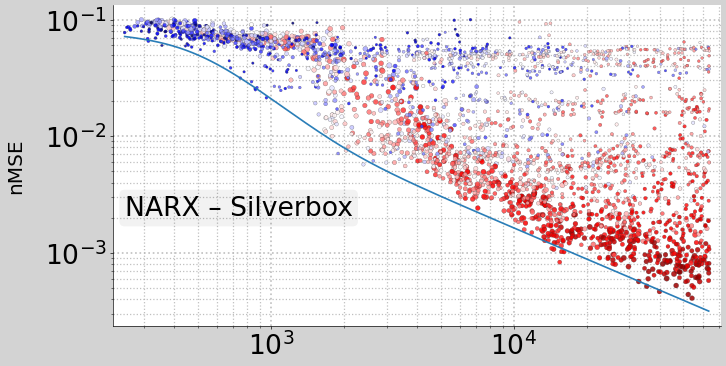}
    \end{minipage}\hfill
    \begin{minipage}[b]{0.48\textwidth}
        \centering
        \includegraphics[width=\linewidth]{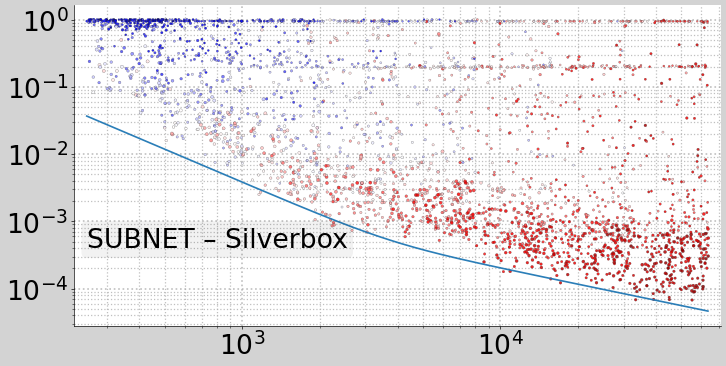}
    \end{minipage}

    \begin{minipage}[b]{0.48\textwidth}
        \centering
        \includegraphics[width=\linewidth]{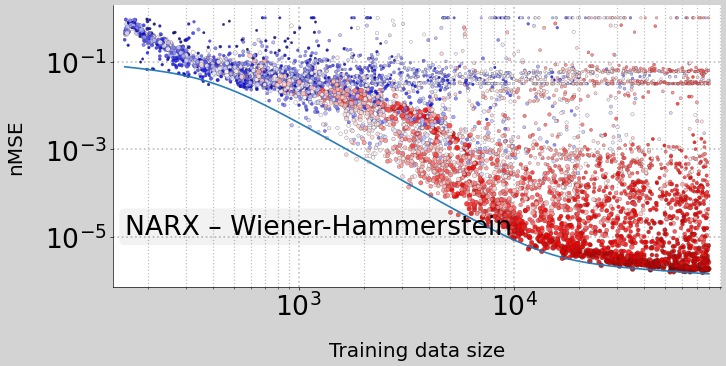}
    \end{minipage}\hfill
    \begin{minipage}[b]{0.48\textwidth}
        \centering
        \includegraphics[width=\linewidth]{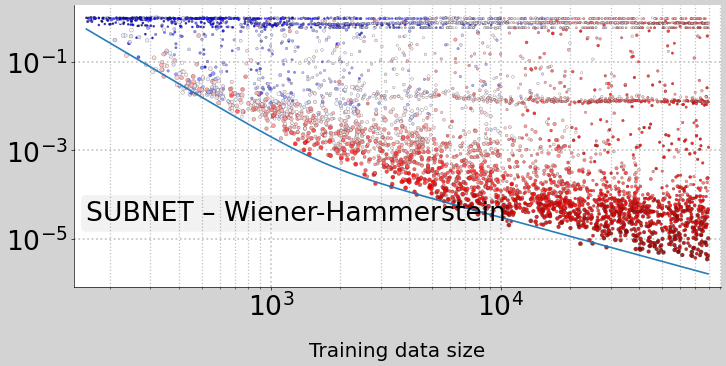}
    \end{minipage}

    \caption{{All data-$\operatorname{nMSE}$ NSLs using input-output data.}}
    \label{fig:data-io-nmse}
\end{figure*}

\begin{figure*}[htbp!] 
    \centering

    \begin{minipage}[b]{0.48\textwidth}
        \centering        \includegraphics[width=\linewidth]{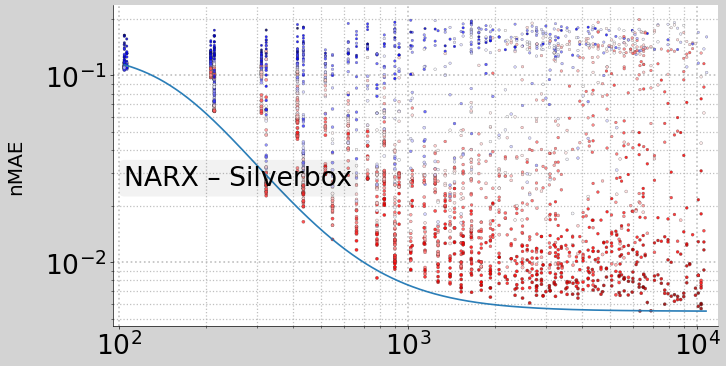}
    \end{minipage}\hfill
    \begin{minipage}[b]{0.48\textwidth}
        \centering
        \includegraphics[width=\linewidth]{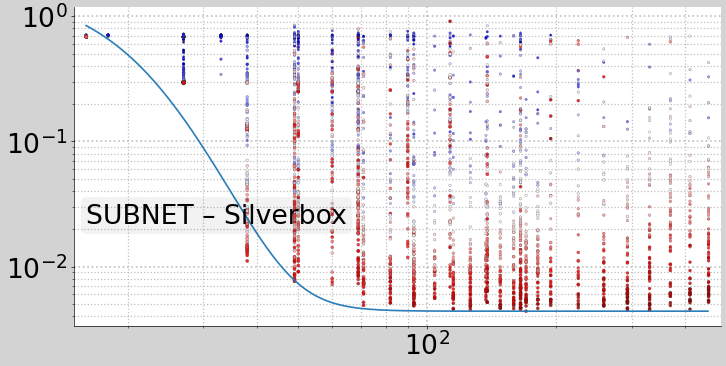}
    \end{minipage}

    \begin{minipage}[b]{0.48\textwidth}
        \centering
        \includegraphics[width=\linewidth]{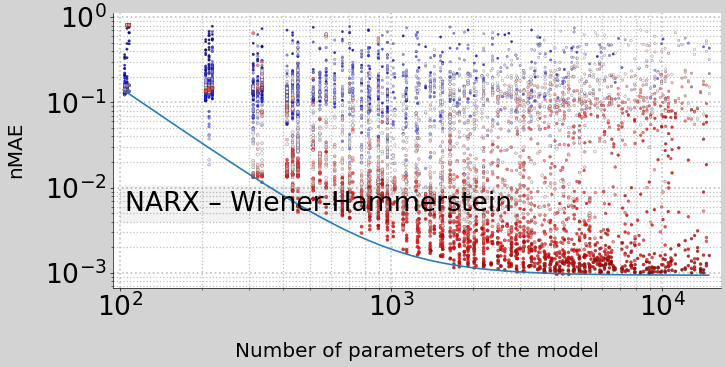}
    \end{minipage}\hfill
    \begin{minipage}[b]{0.48\textwidth}
        \centering
        \includegraphics[width=\linewidth]{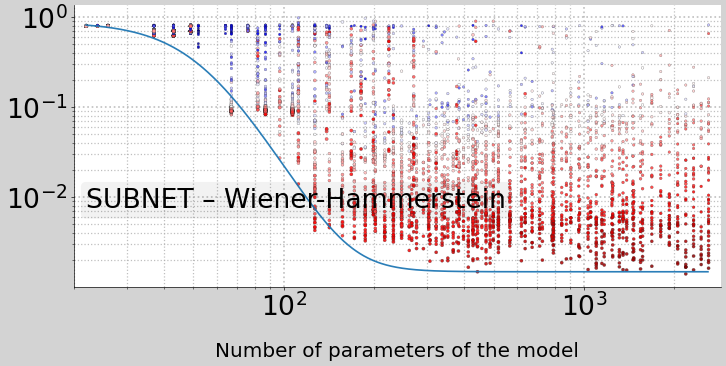}
    \end{minipage}

    \caption{{All model-$\operatorname{nMAE}$ NSLs using input-output data.}}
    \label{fig:model-io-nmae}
\end{figure*}

\begin{figure*}[htbp!] 
    \centering

    \begin{minipage}[b]{0.48\textwidth}
        \centering        \includegraphics[width=\linewidth]{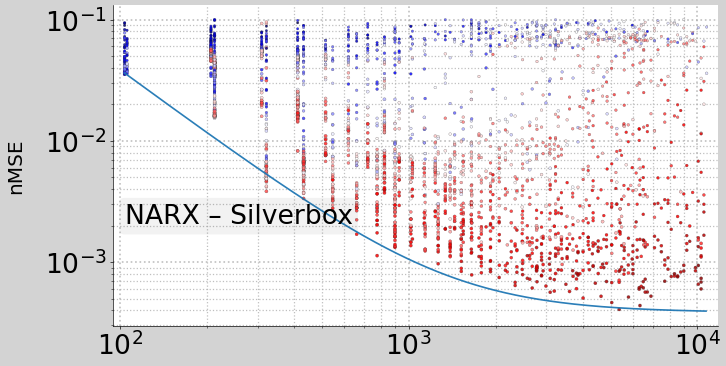}
    \end{minipage}\hfill
    \begin{minipage}[b]{0.48\textwidth}
        \centering
        \includegraphics[width=\linewidth]{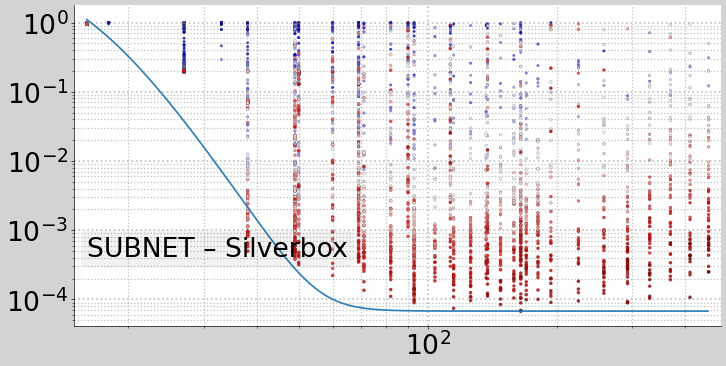}
    \end{minipage}

    \begin{minipage}[b]{0.48\textwidth}
        \centering
        \includegraphics[width=\linewidth]{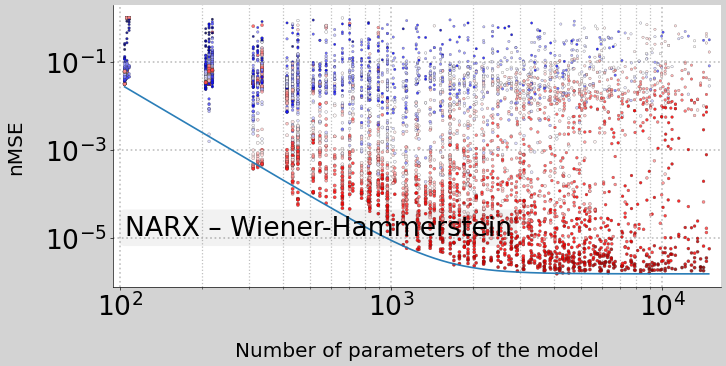}
    \end{minipage}\hfill
    \begin{minipage}[b]{0.48\textwidth}
        \centering
        \includegraphics[width=\linewidth]{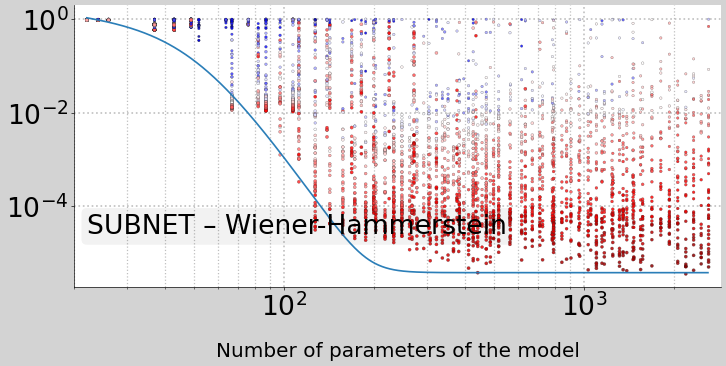}
    \end{minipage}

    \caption{{All model-$\operatorname{nMSE}$ NSLs using input-output data.}}
    \label{fig:model-io-nmse}
\end{figure*}

{In this section, we collect the figures of all considered NSLs using input-output data.
Figures~\ref{fig:comp-io-nmae}~and~\ref{fig:comp-io-nmse} show the compute-nMAE NSLs for input-output data. We can see that they look similar to the corresponding NSLs using input-state-output data in Figures~\ref{fig:comp-nmae}~and~\ref{fig:comp-nmse}. There is one very notable difference, namely that all of these NSLs saturate strongly at more than $10^{12}\, \mathrm{flops}$. On the contrary, there we do not see saturation in the corresponding data-nMAE NSLs in Figures~\ref{fig:data-io-nmae}~and~\ref{fig:data-io-nmse}. This indicates that in order to fully utilize a higher compute than $10^{12}\, \mathrm{flops}$ more data is necessary. We can moreover see that the SUBNET architecture on the Wiener-Hammerstein benchmark has steeper slope than the NARX architecture in the data-NMAE NSLs. This indicates that SUBNET might benefit more from adding more data than NARX on the Wiener-Hammerstein dataset. Inspecting Figures~\ref{fig:model-io-nmae}~and~\ref{fig:model-io-nmse}, we can see that the model NSLs are also highly saturated at relatively small model sizes for input-output data methods.
}

\end{appendices}

\end{document}